\documentclass[12pt]{article}

\usepackage{amsmath,amssymb}
\usepackage{graphicx} 
\makeatletter

\@addtoreset{equation}{section}
\makeatother
\newtheorem{thm}{Theorem}[section]
 \newtheorem{lem}{Lemma}[section]
 \newtheorem{prop}{Proposition}[section]
\begin{document}
\def\be{\begin{equation}}
\def\ba{\begin{eqnarray}}
\def\ee{\end{equation}}
\def\ea{\end{eqnarray}}
\def\la{\langle}
\def\ra{\rangle}
\def\Th{\Theta}
\def\lam{\lambda}
\def\pa{\partial}
\def\R{{\bf R}^2}
\def\oq0{\overline{Q}_0^i}
\def\ve{\varepsilon}
\def\vp{\varphi}
\def\iint{\int\! \! \! \int}
\def\ovw{\overline{w}}
\def\tiw{\tilde{w}}
\def\oxi{\overline{\xi}}
\def\txi{\tilde{\xi}}
\def\tiL{\tilde{\Lambda}_i}
\def\non{\nonumber}
\def\dis{\displaystyle}
\def\var{\varepsilon}
\vspace{1.5cm}
\begin{center}
{\large \bf Lower Bound of the Lifespan of the Solution \\
to Systems of Quasi-linear Wave Equations\\
with Multiple Propagation Speeds}\\[1cm]
\end{center}
\begin{center}
	{\sc Akira HOSHIGA}\\
	Shizuoka University, Japan\\
\end{center}
\vspace{1cm}
       
\begin{abstract}
We consider the Cauchy problem of systems of quasilinear wave equations in 2-dimensional space. We assume that the propagation speeds are distinct and that the nonlinearities  contain quadratic and cubic terms of  the  first and second order derivatives of the solution. We know that if the all quadratic and cubic terms of nonlinearities satisfy $Strong$ $Null$-$condition$, then there exists a global solution for sufficiently small initial data. In this paper, we study about the lifespan of the smooth solution, when the cubic terms in the quasi-linear nonlinearities do not satisfy the Strong null-condition. In the proof of  our claim, we use the $ghost$ $weight$ energy method and the $L^{\infty}$-$L^{\infty}$ estimates of the solution, which is slightly improved. 
\end{abstract}

\noindent
------------------------------------------------------------\\
{\footnotesize 
2010 ASC: 35L05, 35L15, 35L51, 35L70.}

\section{Intrduction}
In this paper, we study the Cauchy problem;
\ba
&\Box_i u^i =\pa_0^2 u^i-c_i^2\triangle u^i=F^i(\pa u, \pa^2 u)\qquad (x, t)\in \mathbb{R}^2\times (0, \infty),
\label{1.1}\\
&u^i(x, 0)=\var f^i(x),\quad \pa_t u^i (x, 0)=\var g^i(x)\qquad\qquad x\in\mathbb{R}^2,\label{1.2}
\ea
where $i=1,2,\cdots, m$ and $u(x, t)=(u^1(x, t), u^2(x, t),\cdots, u^m(x, t))$.  We denote $\pa u=(\pa_\alpha u)_{\alpha=0,1,2}$ with $\pa_0=\pa_t=\pa/\pa t$, $\pa_j=\pa/\pa x_j\ (j=1, 2)$ and $\pa^2u=(\pa_\alpha\pa_{\beta}u)_{\alpha,\beta=0,1,2}$. Let $\var >0$ is a small parameter and assume  $f^i(x), g^i(x)\in C^{\infty}_0(\mathbb{R}^2)$ and $\mbox{supp}\{f^i\}, \mbox{supp}\{g^i\}\subset \{x\in\mathbb{R}^2\ :\ |x|\leq M\}$ for some positive constant $M$. We also assume that the propagation speeds of (\ref{1.1}) are distinct constants, namely we assume
\ba
0<c_1<c_2<\cdots <c_m.\label{1.3}
\ea
Each nonlinearity $F^i(v, w)$ is smooth near the origin and is expressed as 
\ba
F^i(\pa u,\pa^2 u)=\sum_{\ell=1}^m\sum_{\alpha, \beta=0}^2 A_{\ell}^{i, \alpha\beta}(\pa u)\pa_{\alpha}\pa_{\beta}u^{\ell}+B^i(\pa u),\label{1.4}
\ea
where
\ba
 A_{\ell}^{i, \alpha\beta}(\pa u)=\sum_{j=1}^m\sum_{\gamma=0}^2a_{\ell j}^{i, \alpha\beta\gamma}\pa_{\gamma}u^j+\sum_{j,k=1}^m\sum_{\gamma, \delta=0}^2c_{\ell jk}^{i, \alpha\beta\gamma\delta}\pa_{\gamma}u^j\pa_{\delta}u^k+O(|\pa u|^3)\label{1.5}
\ea
and 
\ba
B^i(\pa u)=\sum_{j,k=1}^m\sum_{\alpha,\beta=0}^2b_{jk}^{i, \alpha\beta}\pa_{\alpha}u^j\pa_{\beta}u^k+\sum_{j,k,\ell =1}^m\sum_{\alpha,\beta, \gamma=0}^2d_{jk\ell}^{i, \alpha\beta\gamma}\pa_{\alpha}u^j\pa_{\beta}u^k\pa_{\gamma}u^{\ell}+O(|\pa u|^4).\label{1.6}
\ea
Here $a_{\ell j}^{i, \alpha\beta\gamma}, b_{jk}^{i, \alpha\beta}, c_{\ell jk}^{i, \alpha\beta\gamma\delta}, d_{jk\ell}^{i, \alpha\beta\gamma}$ are constants. \\
In order to derive  energy estimate, we need to assume that for each $i, \ell=1,2,\cdots, m$ and $\alpha,\beta=0,1,2$,
\ba
A_{\ell}^{i,\alpha\beta}(v)=A_i^{\ell,\alpha\beta}(v)=A_i^{\ell,\beta\alpha}(v)\label{1.7}
\ea
and 
\ba
|A_{\ell}^{i,\alpha\beta}(v)|<\frac{(\min\{1, c_1\})^2}{2m}.\label{1.8}
\ea
The assumption (\ref{1.8}) constitutes no additional restriction, since we will only deal with small solutions. Note that by (\ref{1.3}) and (\ref{1.4}), we have 
for any $i=1,2,\cdots,m$,
\ba
u^i(x, t)=0\qquad \mbox{for}\quad |x|\geq c_mt+M.\label{1.85}
\ea
For the proof of (\ref{1.85}), see Theorem 4a in F. John \cite{j1}.
Furthermore, in order to derive  the $ghost$ $weight$ energy method, we need to assume that 
\ba
\begin{aligned}
a_{\ell j}^{i, \alpha\beta\gamma}=0&\qquad\mbox{when}\qquad  (j, \ell) \ne (i, i)\\
b_{jk}^{i, \alpha\beta}=0&\qquad \mbox{when}\qquad j\ne k
\end{aligned}
\label{1.0}
\ea 
for each $i=1,2,\cdots, m$ and $\alpha,\beta,\gamma=0,1,2$. This assumption means that only the terms $\pa u^i\pa^2 u^i$ and $\pa u^j\pa u^j$ $(j=1,\cdots, m)$ appear in the quadratic terms of $F^i$. \\
\indent
Our purpose of this paper is to show a precise estimate for the $lifespan$ $T_{\ve}$. Here, we define $T_{\ve}$ by the supremum of $T>0$ for which there exists a solution 
$u$ to the Cauchy problem (\ref{1.1}) and (\ref{1.2}) in $\big(C^{\infty}(\mathbb{R}^2\times [0, T))\big)^m$. To state the known results and our our result, we introduce some notations. 
Firstly, for $X=(X_0, X_1, X_2)\in \mathbb{R}^3$, we define  $\Phi (X)=(\Phi_{\ell}^i(X))_{i, \ell=1,2,\cdots, m}$, $\Psi (X)=(\Psi_{\ell}^i(X))_{i, \ell=1,2,\cdots, m}$, $\Theta (X)=(\Theta_{\ell}^i(X))_{i, \ell=1,2,\cdots, m}$ and $\Xi (X)=(\Xi_{\ell}^i(X))_{i, \ell=1,2,\cdots, m}$ by
\ba
&&\Phi_{\ell}^i(X)=\sum_{\alpha,\beta,\gamma=0}^2a_{\ell\ell}^{i,\alpha\beta\gamma}X_{\alpha}X_{\beta}X_{\gamma},
\label{1.9}
\\
&&\Psi_{\ell}^i(X)=\sum_{\alpha,\beta=0}^2b_{\ell\ell}^{i,\alpha\beta}X_{\alpha}X_{\beta},
\label{1.10}
\\
&&\Theta_{\ell}^i(X)=\sum_{\alpha,\beta,\gamma,\delta=0}^2c_{\ell\ell\ell}^{i,\alpha\beta\gamma\delta}X_{\alpha}X_{\beta}X_{\gamma}X_{\delta},
\label{1.11}
\\
&&\Xi_{\ell}^i(X)=\sum_{\alpha,\beta,\gamma=0}^2d_{\ell\ell\ell}^{i,\alpha\beta\gamma}X_{\alpha}X_{\beta}X_{\gamma}.
\label{1.12}
\ea
Moreover, let  $\phi (X)=(\phi_{\ell}^i(X))_{i, \ell=1,2,\cdots, m}$ be a function of $X=(X_0, X_1, X_2)$. If 
\ba
\phi^i_{\ell}(X)\equiv 0 \quad\mbox{for}\quad X_0^2=c_{\ell}^2(X_1^2+X_2^2)\label{1.13}
\ea
holds for each $i, \ell=1,2,\cdots,m$, then we denote $\phi\approx 0$ and we say that $\phi$ satisfies $Strong$ $Null$-$condition$. On the other hand, if (\ref{1.13}) holds when $\ell=i$ $(i=1,2,\cdots, m)$, then we denote $\phi\sim 0$ and we say that $\phi$ satisfies $Standard$ $Null$-$condition$. In \cite{2005}, the author showed that 
$\dis{\liminf_{\ve\to 0} \ve \sqrt{T_{\ve}}\geq C}$ holds for a certain positive constant $C$, provided $B^i(\pa u)\equiv 0$. On the other hand, the author also showed in \cite{2006} that  $T_{\ve}=\infty$ for sufficiently small $\ve>0$, provided $\Phi$, $\Psi$, $\Theta$ and $\Xi$ satisfy Strong Null-Condition and $A_{\ell}^{i,\alpha\beta}(\pa u)\equiv 0$ holds for $\ell\ne i$. In this paper, we consider the case that $\Phi$ and $\Psi$ 
satisfy Strong Null-condition and $\Xi$ satisfies Standard Null-condition. Namely, we assume $\Phi\approx 0$, $\Psi\approx 0$ and 
$\Xi\sim 0$. \par
\indent
Secondly, we introduce the Friedlander radiation field $\mathcal{F}^i(\rho,\omega)$. Let $u^i_0(x, t)$ be the solution to the Cauchy problem of the homogeneous linear wave equation; \\
\ba
&\Box_i u_0^i =0\qquad\qquad  (x, t)\in \mathbb{R}^2\times (0, \infty),
\label{1.14}\\
&u_0^i(x, 0)=f^i(x),\quad \pa_0 u_0^i (x, 0)=g^i(x)\qquad x\in\mathbb{R}^2.\label{1.15}
\ea
Then we define $\mathcal{F}^i$ by
\ba
\mathcal{F}^i(\rho, \omega)=\lim_{r\to\infty}r^{\frac12}u_0^i(x, t)\label{1.16}
\ea
with $x=r\omega \ (\omega\in S^1)$ and $\rho=r-c_it$. We know that $\mathcal{F}^i(\rho,\omega)$ is expressed by
\ba
\non
\mathcal{F}^i(\rho, \omega)=\frac1{2\sqrt{2}}\int_{\rho}^{\infty}\frac1{\sqrt{s-\rho}}(R_{g^i}(s, \omega)-\pa_s
R_{f^i}(s, \omega))ds,
\ea
where $R_h(s, \omega)$ is the Radon transform of $h\in C_0^{\infty}(\mathbb{R}^2)$, $i. e.$, 
\ba
\non
R_h(s, \omega)=\int_{\omega\cdot y=s}h(y)\ dS_y
\ea
for $s\in \mathbb{R},\ \omega\in S^1$. Note that $\mathcal{F}^i(\rho, \omega)$ satisfies
\ba
\mathcal{F}^i(\rho, \omega)=0\qquad \mbox{for}\qquad \rho\geq M,\label{1.17}
\ea
\ba
\big|\pa_{\rho}^{\ell}\mathcal{F}^i(\rho,\omega)\big|\leq \frac{C}{(1+|\rho|)^{\frac12+\ell}}\label{1.18}
\ea
and
\ba
\big|r^{\frac12}\pa_0^{\ell}u_0^i(r\omega, t)-(-c_i)^{\ell}\pa_{\rho}^{\ell}\mathcal{F}^i(r-c_it, \omega)\big|\leqq\frac{C(1+|r-c_it|)^{\frac12}}t
\qquad \label{1.19}
\ea
for $t\geq r/(2c_i)\geq 1$ and $\ell=1, 2$. For the details about (\ref{1.17}), (\ref{1.18}) and (\ref{1.19}), see L. H\"ormander \cite{hor1}.\\ 
Then we define a constant
\ba
H_i=\max_{\begin{subarray}{c}\rho\in\mathbb{R}\\
\omega\in S^1\end{subarray}}\left\{-\frac1{c_i^2}\Theta_i^i(-c_i, \omega )
\pa_{\rho}\mathcal{F}^i(\rho,\omega)\pa_{\rho}^2\mathcal{F}^i(\rho,\omega)\right\}\label{1.20}
\ea
and set
\ba
H=\max\{H_1, H_2,\cdots,H_m\}.
\ea
Note that by (\ref{1.17}) and (\ref{1.18}), each $H_i$ is well-defined and nonnegative. \\
\indent
Now, we state our main result.\\
\begin{thm}
Assume that  (\ref{1.7}), (\ref{1.8}) and (\ref{1.0}) hold for the Cauchy problem (\ref{1.1}) and (\ref{1.2}). Also assume $\Phi\approx 0$, $\Psi\approx 0$ and 
$\Xi\sim 0$. Then, if $H>0$, we have
\ba
\liminf_{\varepsilon\to 0}\varepsilon^2\log (1+T_{\varepsilon})\geqq \frac1H.\label{1.21}
\ea
\end{thm}
Note that the  author showed the same estimate when $a_{\ell j}^{i, \alpha\beta,\gamma}=b_{jk}^{i,\alpha\beta}=0$ for all $\alpha, \beta, \gamma=0,1,2$ and $i, j, k=1,2,\cdots, m$ in \cite{2000}. Hence, our result (\ref{1.21}) is a generalization of the result in \cite{2000}. Also note that we can not improve the estimate (\ref{1.21}), in general, since  the counter result has been shown when $m=1$ and $B^1(\pa u)\equiv 0$ in \cite{1995}.
\par
\indent
In the following sections, we aim at showing (\ref{1.21}). In section 2, we prepare some notations and state a lemma which implies (\ref{1.21}).  We also discus about the  estimates of the null-form. In section 3, we will show the $L^{\infty}$-$L^{\infty}$ estimates of solutions to the wave equation. It is an improvement of the one showed in  \cite{hk1}. In section 4, we concentrate to show $a$ $priori$ estimates of the solution, by using the ghost energy inequality  and the method of ordinary differential equation along the characteristic curves.  \\

\section{Preliminary for the proof of Theorem 1.1}
\indent
Our main theorem is immediately derived from the following lemma. 
\begin{lem}
Under the same situation as Theorem 1.1,  
choose a positive constant $B$ to be $B<1/H$. Then there exists a constnat $\ve_0(B)>0$ such that 
\ba
\ve^2\log (1+T_{\ve})\geq B\label{2.0}
\ea
holds for $0<\ve<\ve_0(B)$. 
\end{lem}
In order to state another lemma which causes Lemma 2.1,  we introduce  some notations.  At first, we introduce the following differential operators, 
\ba
\non
\Omega=x_1\pa_2-x_2\pa_1,\qquad S=t\pa_0+x_1\pa_1+x_2\pa_2
\ea
and denote
$$
\Gamma =(\Gamma_0,\ \Gamma_1,\ \Gamma_2,\ \Gamma_3,\ \Gamma_4)=(\pa_0,\ \pa_1,\ \pa_2,\ \Omega,\ S)
$$
and
\ba
 \non
 \Gamma^a=\Gamma_0^{a_0}\Gamma_1^{a_1}\Gamma_2^{a_2}\Gamma_3^{a_3}\Gamma_4^{a_4}
 \ea
 for a multi-index $a=(a_0, a_1,a_2,a_3,a_4)$. 
We can verify the following commutator relations;
\ba
\non
[\Gamma_{\alpha}, \Box_i]=-2\delta_{\alpha 4}\Box_i \qquad (\alpha=0,1,2,3,4,  \ i=1, 2, \cdots, m)
\ea
and
\ba
\non
&[\pa_{\alpha}, \pa_{\beta}]=0\quad (\alpha,\beta=0,1,2),\qquad [S, \pa_{\alpha}]=-\pa_{\alpha}\quad (\alpha=0,1,2)&\\
\non
&[\Omega, \pa_1]=-\pa_2,\quad [\Omega, \pa_2]=\pa_1,\quad [\Omega, \pa_0]=0,\quad [S, \Omega]=0.&
\ea
 Here, $[A, B]=AB-BA$ and $\delta_{\alpha\beta}$ is the Kronecker delta. \\
 \indent
 Secondly, we define norms. Let $v(x, t)=(v^1(x, t), v^2(x, t),\cdots,v^m(x, t))$ be a vector valued function defined on $\mathbb{R}^2\times [0, T)$, then we set
 \ba
 \non
 |v(x, t)|_k&=&\sum_{i=1}^m|v^i(x, t)|_k=\sum_{i=1}^m\sum_{|a|\leq k}|\Gamma^a v^i(x, t)|,\\
  \non
 |v(t)|_k&=&\sup_{x\in \mathbb{R}^2}|v(x, t)|_k,\qquad |v|_{k,T}=\sup_{0\leqq t<T}|v(t)|_k,\\
 \non
  [v(x, t)]_k&=&\sum_{i=1}^{m}[v^i(x, t)]_k=\sum_{i=1}^m\sum_{|a|\leq k}|(1+|x|)^{\frac12}(1+||x|-c_it|)^{\frac{15}{16}}|\Gamma^a v^i(x, t)|,\\
  \non
  [v(t)]_k&=&\sup_{x\in\mathbb{R}^2}[v(x, t)]_k,\qquad [v]_{k, T}=\sup_{0\leqq t<T}[v(t)]_k,\\
  \non
  [[v(x, t)]]_k&=&\sum_{i=1}^{m}[[v^i(x, t)]]_k=\sum_{i=1}^m\sum_{|a|\leq k}|(1+|x|)^{\frac12}(1+||x|-c_it|)|\Gamma^a v^i(x, t)|,\\
  \non
  [[v(t)]]_k&=&\sup_{x\in\mathbb{R}^2}[[v(x, t)]]_k,\qquad [[v]]_{k, T}=\sup_{0\leqq t<T}[[v(t)]]_k,\\
  \non
 \la v(x, t)\ra_k&=&\sum_{i=1}^m\la v^i(x, t)\ra_k=\sum_{i=1}^m\sum_{|a|\leq k}(1+|x|+t)^{\frac{7}{16}}|\Gamma^a v^i(x, t)|,\\
 \non
\la v(t)\ra_k&=& \sup_{x\in \mathbb{R}^2}\la (x, t)\ra_k, \qquad \la v\ra_k=\sup_{0\leq t<T}\la v(t)\ra_k,\\
  \non
 \la\la v(x, t)\ra\ra_k&=&\sum_{i=1}^m\la\la v^i(x, t)\ra\ra_k=\sum_{i=1}^m\sum_{|a|\leq k}(1+|x|+t)^{\frac{1}{2}}|\Gamma^a  v^i(x, t)|,\\
 \non
\la\la v(t)\ra\ra_k&=& \sup_{x\in \mathbb{R}^2}\la\la (x, t)\ra\ra_k, \qquad \la\la v\ra\ra_{k,T}=\sup_{0\leq t<T}\la\la v(t)\ra\ra_k,\\
\non
||v(t)||_k&=&  \sum_{i=1}^{m}\sum_{|a|\leq k}\bigg(\int_{\mathbb{R}^2}|\Gamma^a v^i(x, t)|^2\ dx\bigg)^{\frac12},\qquad ||v||_{k, T}=\sup_{0\leqq t<T}||v(t)||_k,
 \ea
 where $k$ is a nonnegative integer and $|a|=a_0+a_1+\cdots+a_4$ for a multi-index $a=(a_0,a_1,a_2,a_3,a_4)$. \\
 Then, we find that the following lemma implies Lemma 2.1.\\
 \begin{lem}
Let $u(x, t)=(u^1(x, t), u^2(x, t),\cdots, u^m(x, t))\in (C^{\infty}(\mathbb{R}^2\times [0, T)))^m$ be a solution to (\ref{1.1}) and (\ref{1.2}). Choose an integer $k$ so that $k\geq 21$. Let $B>0$ be a constant so that $B<1/H$ and also let $J>0$ be a constant. Then, there exist constants $K=K(B)>0$ and $\var_0=\var_0(J, B)>0$ such that, if 
\ba
 [\pa u]_{k, T}+\la u\ra_{k+1, T}\leq J\var
 \label{2.1}
 \ea
 holds for $0<\var<\var_0$, then 
 \ba
 [\pa u]_{k, T_B}+\la u\ra_{k+1, T_B}\leq K\var
 \label{2.2}
 \ea
 holds for the same $\var$. Here, we have set $T_B=\min \{T, t_B\}$ and $t_B=\exp(B/\var^2)-1$. 
 \end{lem}
 \noindent
 {\bf Proof of Lemma 2.1 providing Lemma 2.2:} We show that Lemma 2.2 implies Lemma 2.1 by contradiction. If the statement of Lemma 2.1 is incorrect, there exists a positive constant $B_0(<1/H)$ such that for any $\ve>0$, there exists $\delta=\delta(\ve)>0$ such that 
\ba
\delta^2\log (1+T_{\delta})\leq B_0 \quad\mbox{and}\quad 0<\delta<\ve.\label{2.00}
\ea
 On the other hand,  by the local existence theorem which was shown in 
A. Majda \cite{ma}, we find that there are positive constants $\ve_1$ and $t_{\ve}$ such that  there exists a smooth solution $u(x, t)\in C^{\infty}(\mathbb{R}^2\times [0, t_{\ve}))$ to (\ref{1.1}) and (\ref{1.2}) for $0<\ve<\ve_1$.   Let $L>0$ be a constant
 satisfying $[\pa u(0)]_k+\la u(0)\ra_{k+1}\leq L\ve$  and set $J_0=2\max\{K(B_0), L\}$, where $K(B_0)$ is the constant determined in Lemma 2.2 with $B=B_0$. Then we can define a positive constant $\tau_{\ve}$ by  
\ba
\non
\tau_{\ve}=\sup\{t>0\ :\ t<T_{\ve}\ \ \mbox{and}\ \ [\pa u(t)]_k+\la u(t)\ra_{k+1}\leq J_0\ve\}\ (\leq T_{\ve})
\ea
for each $\ve\in (0, \ve_*)$. Here we have set $\ve_*=\min\{\ve_0(J_0, B_0) , \ve_1\}$. Note that (\ref{2.1}) holds for $\ve\in (0, \ve_*)$ with $J=J_0$ and $T=\tau_{\ve}$. Moreover, by using (\ref{1.7}) and (\ref{1.8}), we can show $\tau_{{\ve}}<T_{{\ve}}$ for each $\ve\in (0, \ve_*)$. (For the detail, see the proof of Lemma 2.1 in \cite{2005}.)   This means that
\ba
 [\pa u]_{k, \tau_{\ve}}+\la u\ra_{k+1,\tau_{\ve}}= J_0{\ve}\label{2.p}
\ea
holds for $\ve\in (0, \ve_*)$. However, as mentioned above,  there exists a  constant $\delta=\delta (\ve_*)$ such that
(\ref{2.00}) holds. In that case, we find that $T_{B_0}=\min\{\tau_{\delta}, t_{B_0}\}=\tau_{\delta}$ and hence that Lemma 2.2 implies
\ba
 [\pa u]_{k,\tau_{\delta}}+\la u\ra_{k+1,\tau_{\delta}}\leq K(B_0){\delta} \leq \frac{J_0}2{\delta}.
\ea
This contradicts to (\ref{2.p}) and therefore we find that the claim of Lemma 2.1 is correct.\\[2mm]
\indent
In the rest of this paper, we aim at showing Lemma 2.2. For this purpose, we prepare a proposition with respect to the $null$-$form$. Set $\dis{c_*=\min_{1\leq i\leq m}\{c_i-c_{i-1}\}/3}$ with $c_0=0$. We see $c_*>0$ from (\ref{1.3}). Also we set
\ba
\Lambda_i(T)=\{ (x, t)\in \mathbb{R}^2\times [0, T)\ :\ ||x|-c_it|\leq c_*t \}
\label{2.3}
\ea 
and
\ba
\Lambda_0(T)=\mathbb{R}^2\times[0, T)\setminus \dis{\bigcup_{i=1}^m\Lambda_i(T)}.\label{2.4}
\ea
We find that $\Lambda_i(T)\cap\Lambda_j(T)=\emptyset$ holds for any $T>0$, if $i\ne j$ and that there exists a constant $C_1>0$ such that
\ba
\frac1{C_1}(1+|x|+t)\leq 1+||x|-c_jt|\leq C_1(1+|x|+t)\qquad (x, t)\in\Lambda_i(T) \label{2.5}
\ea
holds for any $T>0$, if $i\ne j$. \\
\indent
In order to derive a good decay property from the null-form in $\Lambda_i(T)$, we
introduce the following operators;
\ba
Z=(Z_1^i, Z_2^i),\qquad\quad  Z_{\alpha}^i=c_i\pa_{\alpha}+\frac{x_\alpha}{|x|}\pa_0\qquad (i=1,2,\cdots,m,\ \alpha=1,2).\label{2.6}
\ea
Then we find that
\ba
Z_1^i=\frac{c_it-|x|}{t}\pa_1+\frac{x_1S-x_2\Omega}{|x|t},\qquad 
Z_2^i=\frac{c_it-|x|}{t}\pa_2+\frac{x_2S+x_1\Omega}{|x|t}\label{2.7}
\ea
and hence that
\ba
|Z^iv(x, t)|\leq \frac{||x|-c_it|}{t}|\pa v(x, t)|_0+\frac2t|v(x, t)|_1.\label{2.8}
\ea
Now we have the following.\\
\begin{prop}
Let $T>1$ be a constant and let $k$ be a positive integer. Let $v(x, t)=(v^1(x, t),v^2(x, t),\cdots, v^m(x, t))$ and $w(x, t)=(w^1(x, t), w^2(x, t),\cdots, w^m(x, t))$ be functions belonging to $(C^{\infty}(\mathbb{R}^2\times [0, T)))^{m}$. Assume that  $\Phi\approx 0$, $\Psi\approx 0$, $\Xi\sim 0$ and (\ref{1.0}) hold.  Then, there exists a positive constant $C_2$ independent of $T$ such that
\ba 
\begin{aligned}
&\left|\sum_{\alpha,\beta,\gamma=0}^2a_{ii}^{i,\alpha\beta\gamma}\pa_{\gamma}v^i(x, t)\pa_{\alpha}\pa_{\beta}w^i(x, t)\right|_k\\
\leq \  &C_2\sum_{|b+c|\leq k+1}(|Z^i\Gamma^b v^i(x, t)||\Gamma^c\pa^2 w^i (x, t)|+|\Gamma^b \pa v^i(x, t)||Z^i\Gamma^c\pa w^i (x, t)|),\label{2.9}
\end{aligned}
\ea
\ba
\left|\sum_{\alpha,\beta=0}^2b_{jj}^{i, \alpha\beta}\pa_{\alpha}v^j(x, t)\pa_{\beta}v^j(x, t)\right|_k
\leq C_2\sum_{|b+c|\leq k}|Z^j \Gamma^b v^j(x, t)||\Gamma^c\pa v^j(x, t)|,\label{2.10}
\ea
\ba
\begin{aligned}
&\left|\sum_{\alpha,\beta,\gamma=0}^2d_{iii}^{i, \alpha\beta\gamma}\pa_{\alpha}v^i(x, t)\pa_{\beta}v^i(x, t)\pa_{\gamma}v^i(x, t)\right|_k\qquad\qquad\qquad\qquad\qquad\qquad\qquad \\
\leq \  &C_2\sum_{|b+c+d|\leq k}|Z^i \Gamma^b v^i(x, t)||\Gamma^c \pa v^i(x, t)||\Gamma^d \pa v^i(x, t)|
\label{2.11}
\end{aligned}
\ea
and especially
\ba
\begin{aligned}
&\sum_{|a|\leq k}\left|\sum_{\alpha,\beta,\gamma=0}^2\big\{\Gamma^a(a_{ii}^{i,\alpha\beta\gamma}\pa_{\gamma}v^i(x, t)\pa_{\alpha}\pa_{\beta}v^i(x, t))-a_{ii}^{i,\alpha\beta\gamma}\pa_{\gamma}v^i(x, t)\pa_{\alpha}\pa_{\beta}\Gamma^a v^i(x, t)\big\}\right|\qquad\\
\leq \ &C_2\sum_{|b+c|\leq k+1\atop
|b|,|c|\leq k}|Z^i\Gamma^bv^i(x, t)||\Gamma^c\pa v^i(x, t)|\qquad\qquad
\label{2.12}
\end{aligned}
\ea
hold for $i, j=1,2,\cdots,m$.  Moreover, we find that 
\ba
\begin{aligned}
&\left|\sum_{\alpha,\beta=0}^2b_{jj}^{i, \alpha\beta}\pa_{\alpha}v^j(x, t)\pa_{\beta}v^j(x, t)\right|_k\qquad\qquad\qquad
\\
\leq\  &C_2\bigg(\frac{||x|-c_jt||\pa v^j(x, t)|_{[\frac{k}2]}|\pa v^j(x, t)|_{k}}{1+|x|+t}+\frac{P_{k}(v^j)(x, t)}{1+|x|+t}\bigg)\label{2.13}
\end{aligned}
\ea
holds for $(x, t)\in \Lambda_j(T)\cap\{(y, s) : s\geq 1\}$, $i,j=1,2,\cdots,m$, and that
\ba
\begin{aligned}
&\left|\sum_{\alpha,\beta,\gamma=0}^2a_{ii}^{i, \alpha\beta\gamma}\pa_{\gamma}v^i(x, t)\pa_{\alpha}\pa_{\beta}w^i(x, t)\right|_k\qquad\qquad\qquad
\qquad\qquad\qquad\qquad\qquad\\
\leq\ &C_2\bigg(\frac{||x|-c_it|(|\pa v^i(x, t)|_{[\frac{k}2]}|\pa w^i(x, t)|_{k+1}+|\pa w^i(x, t)|_{[\frac{k+1}2]}|\pa v^i(x, t)|_{k})}{1+|x|+t}+\label{2.14}\\
\ &\qquad +\frac{P_{k}(v^i, w^i)(x, t)}{1+|x|+t}\bigg),
\end{aligned}
\ea
\ba
\begin{aligned}
&\left|\sum_{\alpha,\beta,\gamma=0}^2d_{iii}^{i, \alpha\beta\gamma}\pa_{\alpha}v^i(x, t)\pa_{\beta}v^i(x, t)\pa_{\gamma}v^i(x, t)\right|_k\qquad\qquad\\
\leq \ &C_2\bigg(\frac{||x|-c_it||\pa v^i(x, t)|_{[\frac{k}2]}^2|\pa v^i(x, t)|_{k}}{1+|x|+t}+\frac{|\pa v^i(x, t)|_{[\frac{k}2]}P_{k}(v^i)(x, t)}{1+|x|+t}\bigg),\label{2.15}
\end{aligned}
\ea
and
\ba
\begin{aligned}
&\sum_{|a|\leq k}\left|\sum_{\alpha,\beta,\gamma=0}^2\big\{\Gamma^a(a_{ii}^{i,\alpha\beta\gamma}\pa_{\gamma}v^i(x, t)\pa_{\alpha}\pa_{\beta}v^i(x, t))-a_{ii}^{i,\alpha\beta\gamma}\pa_{\gamma}v^i(x, t)\pa_{\alpha}\pa_{\beta}\Gamma^a v^i(x, t)\big\}\right|\ \ \ \\
\leq\ &C_2\bigg(\frac{||x|-c_it||\pa v^i(x, t)|_{[\frac{k+1}2]}|\pa v^i(x, t)|_{k}}{1+|x|+t}+\frac{Q_{k}(v^i)(x, t)}{1+|x|+t}\bigg)
\qquad
\label{2.16}
\end{aligned}
\ea
hold for $(x, t)\in \Lambda_i(T)\cap\{(y, s) : s\geq 1\}$, $i=1,2,\cdots,m$. Here we have set
\ba
\non
&&P_k(v)(x, t)=|v(x, t)|_{[\frac{k}2]+1}|\pa v(x, t)|_{k}+|\pa v(x, t)|_{[\frac{k}2]}|v(x, t)|_{k+1},\\
\non
&&P_k(v, w)(x, t)=|v(x, t)|_{[\frac{k}2]+1}|\pa w(x, t)|_{k+1}+|\pa v(x, t)|_{[\frac{k}2]}|\pa w(x, t)|_{k+1}+\\
&&\non\qquad\qquad\qquad\qquad+|\pa w(x, t)|_{[\frac{k}2]+1}|\pa v(x, t)|_{k}+|\pa w(x, t)|_{[\frac{k}2]+1}|v(x, t)|_{k+1},\\
\non
&&Q_k(v)(x, t)=|v(x, t)|_{[\frac{k}2]+1}|\pa v(x, t)|_{k}+|\pa v(x, t)|_{[\frac{k}2]}|\pa v(x, t)|_{k}+\\
\non
&&\qquad\qquad\qquad\quad  +|\pa v(x, t)|_{[\frac{k}2]+1}|v(x, t)|_{k+1}.
\ea
\end{prop}
For the proof of Proposition 2.1, see Proposition 2.1 in \cite{2006}.

\section{$L^{\infty}$-$L^{\infty}$ estimate}
In this section, we will show a weighted $L^{\infty}$-$L^{\infty}$ estimate of solutions to inhomogeneous wave equations. It is an improvement of the estimate in Proposition 4.2 in \cite{hk1}. Let $c$ and $T$  be positive constants and $F$ be a function in $C^{\infty}(\mathbb{R}^2\times [0, T))$. Then we introduce an operator $L_{c}(F)$;
\ba
L_c(F)(x, t)=\frac1{2\pi c}\int_0^t\bigg(\int_{|x-y|<c(t-s)}\frac{F(y, s)}{\sqrt{c^2(t-s)^2-|x-y|^2}} dy\bigg)ds
\label{3.1}
\ea
for $(x, t)\in \mathbb{R}^2\times [0, T)$. We know that $L_{c}(F)$ is the solution to the Cauchy problem;
\ba
\non
&(\pa_0^2-c^2\triangle)L_{c}(F)=F(x, t),&\qquad\quad (x, t)\in \mathbb{R}^2\times [0, T),\\
\non
&L_{c}(x, 0)=\pa_0 L_{c}(x, 0)=0,& \qquad\qquad x\in \mathbb{R}^2.
\ea
Then we have the following.
\begin{prop}
Let $c_i$ $(i=1, 2, \cdots, m)$ be the propagation speeds defined in (\ref{1.3}). Let $T>0$ and $F, G ,H \in C^{\infty}(\mathbb{R}^2\times [0, T))$. Choose $\mu>0$, $\nu >0$ and $\rho>0$ arbitrarily. Then, there exist positive constants $\tilde{C}_{\mu}$, $\hat{C}_{\nu}$ and $\dot{C}_{\rho}$ independent of $T$ such that
\ba
|L_{c_i}(F)(x, t)|(1+|x|)^{\frac12}\leq \tilde{C}_{\mu}M^{(i)}_{\mu, 0}(F)(x, t)
\label{l0}
\ea
and
\ba
\begin{aligned}
&|\nabla L_{c_i}(G+H)(x, t)|(1+|x|)^{\frac12}(1+||x|-c_it|)\\
\leq\ \ &\hat{C}_{\nu}M^{(i)}_{\nu, 0}(G)(x, t)+\dot{C}_{\rho}\{1+\log (1+t+|x|)\}M^{(i)}_{0, \rho}(H)(x, t)
\label{l1}
\end{aligned}
\ea
hold for $(x, t)\in \R\times [0, T)$. Here $\nabla=(\pa_1, \pa_2)$ and we have set
\ba
\non
M^{(i)}_{\mu,\nu}(F)(x, t)=\sum_{j=0}^m\sup_{(y, s)\in \atop \Lambda_j(T)\cap D^i(x, t)}\{|y|^{\frac12}z_{\mu, \nu}^{(j)}(|y|, s)|F(y, s)|_1\},
\label{l2}
\ea
\ba
\non
z_{\mu, \nu}^{(j)}(\lam, s)=(1+s+\lam)^{1+\mu}(1+|\lam -c_j s|)^{1+\nu},
\label{l3}
\ea
\ba
\non
D^i(x, t)=\{(y, s)\ :\ |x-y|\leq  c_i(t-s)\ \}.
\ea
\end{prop}

\noindent
{\bf Proof of Proposition 3.1:}\ \ \ \ By the same argument with the proof of Propositions 4.1 and 4.2 in \cite{hk1}, we obtain (\ref{l0}) and (\ref{l1}) when $H(x, t)\equiv 0$. Therefore, we have only to show (\ref{l1}) when $G(x, t)\equiv 0$. Without loss of generality, we may assume $c_i=1$ and for the sake of simplicity, we denote the constant depending on $\rho$ by $C$ which may change line by line,
during this section.  Set 
\ba
\non
E_1&=&\{(y, s)\in \R\times[0, t)\ :\ |y|+s>t-r,\  |x-y|<t-s\}\\
\non
E_2&=&\{(y, s)\in \R\times[0, t)\ :\ (t-r-1/2)_+<|y|+s<t-r \}\\
\non
E_3&=&\{(y, s)\in \R\times[0, t)\ :\ |y|+s<(t-r-1/2)_+ \}
\ea
with $r=|x|$ and define
\ba
\non
P_j(H)(x, t)=\frac1{2\pi}\iint_{E_j}\frac{H(y, s)}{\sqrt{(t-s)^2-|x-y|^2}}\ dyds\qquad \qquad (j=1,2,3),
\ea
then we have
$$
\pa_{\ell} L_{c_i}(H)(x, t)=\sum_{j=1}^3P_j(\pa_{\ell} H)(x, t)\qquad \quad (\ell=1, 2.)
$$\\
\indent
Firstly, we deal with $P_1(\pa_{\ell} H)(x, t)$. Following the computation made in Section 4 of \cite{hk1}, we  find 
\ba
\non
|P_1(\pa_{\ell} H)(x, t)|\leq CM^{(i)}_{0,\rho}(H)\sum_{k=1}^5I_k,
\label{l4}
\ea
where we have set
\ba
\non
I_1&=&\sum_{j=0}^m \iint_{D_1}\frac{\lam^{\frac12}}{z_{0, \rho}^{(j)}(\lam, s)}\left(\int_{-\varphi}^{\varphi}K_1(\lam, \psi;
r, t-s)\  d\psi\right)d\lam ds, \\
\non
I_2&=&\sum_{j=0}^m \int_{D'_2}\frac{\lam^{\frac12}}{z_{0, \rho}^{(j)}(\lam, s)}\left(\int_{0}^{1}K_2(\lam, \tau;
r, t-s)\  d\tau\right)d\sigma, \\
\non
I_3&=&\sum_{j=0}^m \iint_{D_2}\frac1{\lam^{\frac12}z_{0, \rho}^{(j)}(\lam, s)}\left(\int_{0}^{1}K_2(\lam, \tau;
r, t-s)\  d\tau\right)d\lam ds, \\
\non
I_4&=&\sum_{j=0}^m \iint_{D_2}\frac{\lam^{\frac12}}{z_{0, \rho}^{(j)}(\lam, s)}\left(\int_{0}^{1}|\pa_{\lam}K_2(\lam, \tau;
r, t-s)|\  d\tau\right)d\lam ds, \\
\non
I_5&=&\sum_{j=0}^m \iint_{D_2}\frac{\lam^{\frac12}}{z_{0, \rho}^{(j)}(\lam, s)}\left(\int_{0}^{1}|(\pa_{\lam}\Psi\cdot K_2)(\lam, \tau;
r, t-s)|\  d\tau\right)d\lam ds.
\ea
Here we have used the following notation:
\ba
\non
K_1(\lam, \varphi; r,t)&=&\frac1{2\pi\sqrt{t^2-r^2-\lam^2+2r\lam\cos \psi}},\\
\non
K_2(\lam, \tau; r,t)&=&\frac1{2\pi\sqrt{2r\lam\tau (1-\tau)\{2-(1-\cos \varphi)\tau\}}},\\
\non
\varphi(\lam; r,t)&=&\arccos \left(\frac{r^2+\lam^2-t^2}{2r\lam}\right),\\
\non
\Psi(\lam, \tau; r,t)&=&\arccos \{1-(1-\cos\varphi)\tau\},\\
\non
D_1&=&\{(\lam, s)\in (0, \infty)\times (0, t)\ :\ \lam_-<\lam<\lam_-+\delta \ \ \  \mbox{or}\ \ \ \lam_+-\delta<\lam<\lam_+\},\\
\non
D_2&=&\{(\lam, s)\in (0, \infty)\times (0, t)\ :\ \lam_-+\delta<\lam<\lam_+-\delta \},\\
\non
D'_2&=&\{(\lam, s)\in (0, \infty)\times (0, t)\ :\ \lam=\lam_-+\delta\ \ \ \mbox{or}\ \ \ \lam=\lam_+-\delta \}\\
\non
\ea
with $\lam_-=|t-s-r|$, $\lam_+=t-s+r$ and $\delta=\min\{r,1/2\}$. Thus we aim to show
\ba
I_k\leq\frac{C\{1+\log(1+t+r)\}}{(1+r)^{\frac12}(1+|t-r|)}\qquad\quad  (k=1,2,3,4,5).\label{l5}
\ea
In oder to show (\ref{l5}), we use the following estimates which are proved in Lemma 4.1 in \cite{hk2}.\\
\begin{lem} \ \ \ Let $(\lam, s)\in D_1\cup D_2$. then we have
\ba
\int_{-\varphi}^{\varphi}K_1\ d\psi=2\int_0^1K_2\ d\tau &\leq&\frac{C}{(r\lam)^{\frac12}}\log\left(2+\frac{r\lam h(t-s-r)}{(\lam-\lam_-)(\lam_+-\lam)}\right),\label{l6}\\
\int_0^1|\pa_{\lam}K_2|\ d\tau &\leq&\frac{C}{(r\lam)^{\frac12}(\lam+s+r-t)},\label{l7}\\
\int_0^1|\pa_{\lam}\Psi\cdot K_2|\ d\tau &\leq&\frac{C}{(r\lam)^{\frac12}}\left(\frac1{(\lam_+-\lam)^{\frac12
}(\lam-\lam_-)^{\frac12}}+\frac1{(\lam^2-\lam_-^2)^{\frac12}}\right),\label{l8}
\ea
where, $h(p)=1$ for $p>0$ and $h(p)=0$ for $p\leq 0$.
\end{lem}
First we evaluate $I_1$. When $t-r-s>0$ and $\lam>\lam_+-\delta$, we have 
\ba
\non
\log\left(2+\frac{r\lam}{(\lam-\lam_-)(\lam_++\lam)}\right)\leq \log 3,
\ea
since $\lam-\lam_->r$. Moreover, we find that
\ba
\non
&z_{0, \rho}^{(j)}(\lam, s)\geq z_{0, \rho}^{(j)}(\lam_+,s)&\qquad \quad \mbox{for}\qquad \lam_+-\delta<\lam<\lam_+\\
\non
&z_{0, \rho}^{(j)}(\lam, s)\geq C z_{0, \rho}^{(j)}(\lam_-,s)&\qquad \quad\mbox{for}\qquad \lam_-<\lam<\lam_-+\delta.
\ea
Hence, by (\ref{l6}), we  get
\ba
I_1\leq \frac{C}{r^{\frac12}}\sum_{j=0}^m(A_{1,j}+B_{1, j}+C_{1, j}),\label{l9}
\ea
where we have set
\ba
\non
A_{1,j}&=&\int_0^t\left(\int_{\lam_+-\delta}^{\lam_+}\frac1{z_{0, \rho}^{(j)}(\lam_+, s)}\ d\lam\right)ds,\\
\non
B_{1,j}&=&\int_0^{(t-r)_+}\left(\int^{\lam_-+\delta}_{\lam_-}\frac1{z_{0, \rho}^{(j)}(\lam_-, s)}\log\left(2+\frac{\lam}{\lam-\lam_-}\right)d\lam\right)ds,\\
\non
C_{1,j}&=&\int_{(t-r)_+}^t\left(\int^{\lam_-+\delta}_{\lam_-}\frac1{z_{0, \rho}^{(j)}(\lam_-, s)}\ d\lam\right)ds.
\ea
It follows that 
\ba
\non
A_{1, j}&=&\int_0^t\left(\int_{\lam_+-\delta}^{\lam_+}\frac1{(1+s+\lam_+)(1+|\lam_+-c_js|)^{1+\rho}}\ d\lam\right)ds\\
&\leq& \frac{C\delta}{1+t+r}\int_{-\infty}^{\infty}\frac1{(1+|(1+c_j)s-t-r|)^{1+\rho}}\ ds\label{l10}\\
\non
&\leq& \frac{C\delta}{1+|t-r|}.
\ea
When we deal with $B_{1,j}$, we may assume $t>r$, since $B_{1, j}=0$ if $t\leq r$.  Integrating by parts, we find
\ba
\non
\int^{\lam_-+\delta}_{\lam_-}\log\left(2+\frac{\lam}{\lam-\lam_-}\right)d\lam
&=&
\non
\int^{\lam_-+\delta}_{\lam_-}\left\{\log(3\lam-2\lam_-)-\log(\lam-\lam_-)\right\}d\lam\\
\non
&=&\left[\frac{3\lam-2\lam_-}3\log (3\lam-2\lam_-)-(\lam-\lam_-)\log (\lam-\lam_-)\right]_{\lam_-}^{\lam_-+\delta}\\
\non
&=&\frac{\lam_-+3\delta}{3}\log(\lam_-+3\delta)-\delta\log\delta-\frac{\lam_-}{3}\log\lam_-\\
\non
&=&\frac{\lam_-}3\log \left(1+\frac{3\delta}{\lam_-}\right)+\delta\log (\lam_-+3\delta)-\delta\log{\delta}\\
\non
&\leq& \delta+\delta\log (2+|t-r|)+\delta^{\frac12},
\ea
where we have used $0<\delta<1/2$ and the facts
\ba
\non
0\leq &\dis{\frac{\log (1+x)}{x}}&< 1\qquad\quad \mbox{for}\qquad x> 0,\\
\non
0\leq &-\delta^{\frac12}\log \delta&<1\qquad\quad \mbox{for}\qquad 0<\delta<\frac12.
\ea
Hence we have
\ba
\begin{aligned}
B_{1, j}\ \ \leq& \ \ \frac{C\delta^{\frac12}\log (2+|t-r|)}{1+|t-r|}
\int_{-\infty}^{\infty}\frac1{(1+|(1+c_j)s-t-r|)^{1+\rho}}\ ds\\
\leq&\ \  \frac{C\delta^{\frac12}\{1+\log(1+t+r)\}}{1+|t-r|}.\label{l11}
\end{aligned}
\ea
When $s>(t-r)_+$, we have
\ba
\non
s+\lam_-=2s-t+r\geq |t-r|,\qquad s+\lam_-\geq C|\lam-c_j s|=C|(1-c_j)s-t+r|,
\ea
which imply
\ba
\non
z_{0, \rho}^{(j)}(\lam_-, s)&\geq& C(1+|t-r|)(1+|(1-c_j)s-t+r|)^{1+\rho}\qquad \mbox{if}\quad j\ne i\\
\non
z_{0, \rho}^{(i)}(\lam_-, s)&\geq& C(1+|t-r|)^{1+\rho}(1+2s-t+r).
\ea
Therefore, we get
\ba
\begin{aligned}
C_{1,j}\ \ \leq &\ \ \frac{C\delta}{1+|t-r|}\int_{-\infty}^{\infty} \frac{1}{(1+|(1-c_j)s-t+r|)^{1+\rho}}\ ds\\
\leq\ \ &\frac{C\delta}{1+|t-r|}\qquad \qquad \qquad \mbox{if}\qquad j\ne i
\label{l12}
\end{aligned}
\ea
\ba
\begin{aligned}
C_{1,i}\ \ \leq&\ \ \frac{C\delta}{(1+|t-r|)^{1+\rho}}\int_{(t-r)_+}^{t}\frac1{1+2s+t-r}\ ds\\
\leq&\ \ \frac{C\delta\{1+\log(1+t+r)\}}{1+|t-r|}.\label{l13}
\end{aligned}
\ea
Summing up (\ref{l9}), (\ref{l10}), (\ref{l11}), (\ref{l12}) and (\ref{l13}), we obtain (\ref{l5}) for $k=1$, since $\delta /r\leq 2/(1+r)$. \\
\indent
In the remainder of the proof of (\ref{l4}), we assume $r\geq 1/2$ so that $\delta=1/2$, because $D_2$ is the empty set, if $0<r<1/2$.\\
\indent
 Since $\lam=\lam_-+1/2$ or $\lam=\lam_+-1/2$ for $(\lam, s)\in D_2'$, we obtain (\ref{l5}) for $k=2$ analogously to the previous argument.\\
\indent
Next we evaluate $I_3$. Note that $\lam >\lam_-+1/2$ for $(\lam, s)\in D_2$ and that 
\ba
\non
\log\left(2+\frac{r\lam}{(\lam-\lam_-)(\lam_++\lam)}\right)\leq \log (2+2\lam)
\ea
for $\lam>\lam_-+1/2$. Therefore we get from (\ref{l6})
\ba
\non
r^{\frac12}I_3&\leq& C\sum_{j=0}^m\iint_{D_2}\frac{\log(2+2\lam)}{(1+\lam)z_{0, \rho}^{(j)}(\lam, s)}d\lam ds\\
\non
&\leq& C\{1+\log (1+t+r)\}\sum_{j=0}^m A_{3, j},
\ea
where we have set
\ba
\non
A_{3, j}&=&\iint_{D_2}\frac{1}{(1+s+\lam)^2(1+|\lam-c_js|)^{1+\rho}} d\lam ds\qquad (1\leq j\leq m),\\
\non
A_{3, 0}&=&\iint_{D_2}\frac{1}{(1+s+\lam)(1+\lam)^{2+\rho}} d\lam ds.
\ea
When $1\leq j\leq m$, changing variables by
\ba
\alpha=\lam+s\qquad \mbox{and}\qquad \beta=\lam-s,\label{l14}
\ea
we have
\ba
\non
A_{3, j}&\leq&\frac12\int_{|t-r|}^{t+r}\frac1{(1+\alpha)^2}\left(\int_{-|t-r|}^{\alpha}\frac1{(1+|\psi_j(\alpha,\beta)|)^{1+\rho}}\ d\beta\right)d\alpha\\
\non
&\leq&\frac{C}{1+|t-r|},
\ea
where 
\ba
\non
2\psi_j(\alpha, \beta)=(c_j+1)\beta-(c_j-1)\alpha.
\ea
On the other hand, when $j=0$, we have
\ba
\non
A_{3, 0}\leq \frac{1}{1+|t-r|}\iint_{D_2}\frac1{(1+\lam)^{2+\rho}}\ d\lam ds\leq \frac{C}{1+|t-r|}.
\ea
Therefore we have (\ref{l5}) for $k=3$. \\
\indent
Next we evaluate $I_4$. Since $\lam+s+r-t\geq 1/2$ for $\lam\geq\lam_-+1/2$, we get from (\ref{l7})
\ba
\non
r^{\frac12}I_4&\leq&C\sum_{j=0}^m\iint_{D_2}\frac1{z_{0, \rho}^{(j)}(\lam, s)(\lam+s+r-t+1)}\ d\lam ds\\
\non
&\leq&\frac{C}{(1+|t-r|)}\int_{|t-r|}^{t+r}\frac1{\alpha+r-t+1}\left(\int_{-|r-t|}^{\alpha}\frac1{(1+|\psi_j(\alpha,\beta)|)^{1+\rho}}\ d\beta\right)d\alpha\\
\non
&\leq& \frac{C\{1+\log (1+t+r)\}}{1+|t-r|},
\ea
which yields (\ref{l5}) for $k=4$. \\
\indent
Next we evaluate $I_5$. It follows from $\lam_-+1/2\leq \lam\leq \lam_+-1/2$ that
\ba
\non
3(\lam_+-\lam)\geq \lam_+-\lam+1,\quad
3(\lam-\lam_-)\geq\lam-\lam_-+1,\quad
9(\lam^2-\lam_-^2)\geq (\lam+1)^2-\lam_-^2.
\ea
Hence we get from (\ref{l8})
\ba
\non
r^{\frac12}I_5\leq C\sum_{j=0}^m(A_{5, j}+B_{5, j}+C_{5, j}),
\ea
where we have set
\ba
\non
A_{5, j}&=&\iint_{D_2\cap \{t-r\leq s\}}\frac1{z_{0, \rho}^{(j)}(\lam,s)(t+r-s-\lam)^{\frac12}(\lam-t+s+r)^{\frac12}}
\ d\lam ds,\\
\non
B_{5, j}&=&\iint_{D_2\cap \{t-r\geq s\}}\frac1{z_{0, \rho}^{(j)}(\lam,s)(t+r-s-\lam+1)^{\frac12}(\lam+t-s-r+1)^{\frac12}}
\ d\lam ds,\\
\non
C_{5, j}&=&\iint_{D_2}\frac1{z_{0, \rho}^{(j)}(\lam,s)(\lam-t+s+r+1)^{\frac12}(\lam+t-s-r+1)^{\frac12}}
\ d\lam ds.
\ea
Changing variables by (\ref{l14}), we have
\ba
\non
A_{5, j}&\leq&\frac{C}{1+|t-r|}\int_{|t-r|}^{t+r}\frac1{(t+r-\alpha)^{\frac12}(\alpha-t+r)^{\frac12}}
\left(\int_{-|r-t|}^{\alpha}\frac1{(1+|\psi_j(\alpha,\beta)|)^{1+\rho}} \ d\beta\right)d\alpha\\
\non
&\leq&\frac{C}{1+|t-r|}\int_{t-r}^{t+r}\frac1{(t+r-\alpha)^{\frac12}(\alpha-t+r)^{\frac12}}\ d\alpha\\
\non
&\leq&\frac{C}{1+|t-r|}.
\ea
Changing variables by (\ref{l14}) and then by $\sigma=\psi_j(\alpha, \beta)$, we get
\ba
\non
B_{5, j}
&\leq&\frac12\int_{|t-r|}^{t+r}\frac1{(1+\alpha)(t+r-\alpha+1)^{\frac12}}\times\\
\non
&&\qquad\qquad\qquad\times\left(\int_{\gamma_j}^{\alpha}\frac1{(1+|\sigma|)^{1+\rho}\{1+\frac2{c_j+1}(\sigma-\gamma_j)\}^{\frac12}}\ d\sigma\right)d \alpha,
\ea
where 
\ba
\non
2\gamma_j=(1-c_j)\alpha+(1+c_j)(r-t).
\ea
It has been shown in Lemma 3.13 in \cite{ky} that
\ba
\non
\int_{\gamma_j}^{\alpha}\frac1{(1+|\sigma|)^{1+\rho}\{1+\frac2{c_j+1}(\sigma-\gamma_j)\}^{\frac12}}\ d\sigma
\leq \frac{C}{(1+|\gamma_j|)^{\frac12}}.
\ea
Therefore, if $j\ne i$, we have
\ba
\non
B_{5, j}&\leq&\frac{C}{(1+|t-r|)}\int_{|t-r|}^{t+r}\frac1{(t+r-\alpha+1)^{\frac12}(1+|\gamma_j|)^{\frac12}}\ d\alpha\\
\non
&\leq&\frac{C}{(1+|t-r|)}\int_{|t-r|}^{t+r}\left(\frac1{t+r-\alpha+1}+\frac1{1+|\gamma_j|}\right)\ d\alpha\\
\non
&\leq&\frac{C\{1+\log (1+t+r)\}}{1+|t-r|}.
\ea
On the other hand, if $j=i$, since $\gamma_i=r-t$, we have
\ba
\non
B_{5, i}&\leq&\frac{C}{(1+|t-r|)}\int_{|t-r|}^{t+r}\frac1{(1+\alpha)^{\frac12}(t+r-\alpha+1)^{\frac12}}\ d\alpha\\
\non
&\leq&\frac{C}{(1+|t-r|)}\int_{|t-r|}^{t+r}\left(\frac1{1+\alpha}+\frac1{t+r-\alpha+1}\right)\ d\alpha\\
\non
&\leq&\frac{C\{1+\log (1+t+r)\}}{1+|t-r|}.
\ea
Since we can deal with $C_{5, j}$ similarly to $B_{5, j}$, we obtain (\ref{l5}) for $k=5$. \\
\indent
Secondly, we deal with $P_2(\pa_{\ell} H)$ We can assume $t>r$, since otherwise $E_2$ is empty. Switching to polar coordinates, 
\ba
x=(r\cos \theta, r\sin \theta),\qquad y=\lambda\xi=(\lam \cos (\theta+\psi), \lam\sin (\theta+\psi)),\label{l15}
\ea 
we get
\ba
\non
P_2(\pa_{\ell} H)(x, t)=\int_0^{t-r}\left(\int_{(\lam-\frac12)_+}^{\lam_-}\lam\pa_{\ell} H(\lam\xi, s)\left(\int_{-\pi}^{\pi}K_1(\lam,\psi;r, t-s)\ d\psi\right)d\lam\right)ds
\ea
By Proposition 5.2 in \cite{ay}, we have
\ba
\int_{-\pi}^{\pi}K_1(\lam,\psi;r, t-s)\ d\psi\leq \frac{C}{(\lam+\lam_-)^{\frac12}(\lam_+-\lam)^{\frac12}}\log\left(2+\frac{r\lam}{(\lam_--\lam)(\lam_++\lam)}\right)
\label{l16}
\ea
for $0<s<t-r$ and $0<\lam<\lam_-$. It follows from the fact
\ba
\frac1{\lam_+-\lam}\leq \frac{2}{(r+1)(\lam_--\lam)}\qquad\quad\mbox{for}\qquad 0<s<t-r, \ \ \ \lam_--\frac12\leq\lam<\lam_-\label{l17}
\ea
that
\ba
(r+1)^{\frac12}|P_2(\pa_{\ell} H)(x, t)|\leq CM^{(i)}_{0, \rho}(F)\sum_{j=0}^m A_{6, j},\label{l18}
\ea
where we have set
\ba
\non
A_{6, j}=\int_0^{t-r}\left(\int_{(t-r-\frac12)_+}^{t-r}\frac{\log\left(2+\dis{\frac{\lam+s}{\lam_--\lam}}\right)}{(1+s+\lam)(1+|\lam-c_js|)^{1+\rho}(\lam_-\lam)^{\frac12}}\ d\lam\right)ds.
\ea
Changing variables by (\ref{l14}), we get
\ba
\non
A_{6, j}&\leq&C\int_{t-r-\frac12}^{t-r}\frac{\log\left(2+\dis{\frac{\alpha}{t-r+\alpha}}\right)}{(1+\alpha)
(t-r-\alpha)^{\frac12}}\left(\int_{-\infty}^{\infty}\frac1{(1+|\psi_j|)^{1+\rho}}\ d\beta\right) d\alpha\\
\non
&\leq& \frac{C}{t-r+1/2}\int_{t-r-\frac12}^{t-r}\frac{\log(2(t-r)-\alpha)-\log (t-r-\alpha)}{(t-r-\alpha)^{\frac12}}\ d\alpha.
\ea
Here
\ba
\non
&&\int_{t-r-\frac12}^{t-r}\frac{\log(2(t-r)-\alpha)-\log (t-r-\alpha)}{(t-r-\alpha)^{\frac12}}\ d\alpha\\
\non
&=&\bigg[-2(t-r-\alpha)(\log(2(t-r)-\alpha)-\log(t-r-\alpha))\bigg]_{t-r-\frac12}^{t-r}\\
\non
&&\qquad \qquad -2\int_{t-r-\frac12}^{t-r}\left(\frac{(t-r-\alpha)^{\frac12}}{2(t-r)-\alpha}-\frac{(t-r-\alpha)^{\frac12}}{t-r-\alpha}\right) d\alpha\\
\non
&=&\sqrt{2}\log\bigg(\frac12+t-r\bigg)+\int_{t-r-\frac12}^{t-r}\frac{t-r}{(2(t-r)-\alpha)(t-r-\alpha)^{\frac12}}\ d\alpha\\
\non
&\leq& C\log\bigg(\frac32+2t\bigg)+\frac{t-r}{t-r+\frac12}\int_{t-r-\frac12}^{t-r}\frac{1}{(t-r-\alpha)^{\frac12}}\ d\alpha\\
\non
&\leq& C\bigg(1+\log\bigg(\frac32+2t\bigg)\bigg)
\ea
implies
\ba
A_{6, j}\leq\frac{C\{1+\log (1+t+r)\}}{1+|t-r|}.\label{l19}
\ea
Thus, combining (\ref{l18}) and (\ref{l19}), we obtain 
\ba
\non
(1+r)^{\frac12}(1+|t-r|)|P_2(\pa_{\ell} H)(x, t)|\leq C\{1+\log (1+t+r)\}M^{(i)}_{0, \rho}(H).
\ea
\indent
Thirdly, we deal with $P_3(\pa_{\ell} H)$. We can assume $t>r+1/2$, since otherwise $E_3$ is empty. Integrating by parts in $y$ and switching to polar coordinates as in (\ref{l16}), we get
\ba
\non
P_3(\pa_{\ell} H)(x, t)&=&\int_{0}^{t-r-\frac12}\left(\int_0^{t-r-\frac12}\left(\int_{-\pi}^{\pi}\lam H(\lam\xi, s)K_3(\lam,\psi;x, t-s)d\psi\right)d\lam\right)ds\\
\non
&&+\int_0^{t-r-\frac12}\left(\int_{-\pi}^{\pi}\lam \xi_l H(\lam\xi, s)K_1(\lam, \psi; x, t-s)\bigg|_{\lam=t-s-r-\frac12}d\psi\right)ds\\
\non
&\equiv&J_1+J_2,
\ea
where we have set
\ba
\non
K_3(\lam, \psi; x, t)=\frac{-(x_{\ell}-\lam\xi_{\ell})}{2\pi (t^2-r^2-\lam^2+2r\lam \cos\psi)^{\frac32}}.
\ea
We see from (\ref{l16}) and (\ref{l17}) that
\ba
\non
J_2&\leq& \frac{CM^{(i)}_{0, \rho}(H)\{1+\log (1+t+r)\}}{(r+1)^{\frac12}(t-r+1/2)}\sum_{j=0}^m\int_0^{t-r-\frac12}\frac{1}{(1+|t-r-(c_j+1)s-1/2|)^{1+\rho}}\ ds\\
\non
&\leq&\frac{C\{1+\log (1+t+r)\}M^{(i)}_{0, \rho}(H)}{(1+r)^{\frac12}(1+|t-r|)}.
\ea
As shown in the proof of Proposition 4.2 in \cite{hk2}, we have
\ba
\non
\int_{-\pi}^{\pi}|K_3(\lam, \psi;r, t-s)|d\psi\leq \frac{C}{(\lam_--\lam)(\lam_-+\lam)^{\frac12}(\lam_+-\lam)^{\frac12}}.
\ea
Therefore, since $r+1\leq 2(\lam_+-\lam)$ for $\lam<t-r-s-1/2$, we get from (\ref{l14})
\ba
\non
&&J_1\\
\non
&\leq&\frac{CM^{(i)}_{0, \rho}(H)}{(1+r)^{\frac12}}\sum_{j=0}^m\int_{0}^{t-r-\frac12}\left(\int_0^{t-r-s-\frac12}\frac1{(1+s+\lam)(\lam_--\lam)(1+|\lam-c_js|)^{1+\rho}}\ d\lam\right)ds\\
\non
&\leq&\frac{CM^{(i)}_{0, \rho}(H)}{(1+r)^{\frac12}}\sum_{j=0}^m\left\{\int_{\frac{t-r}2}^{t-r-\frac12}\left(\int_{-\alpha}^{\alpha}\frac1{(1+\alpha)(t-r-\alpha)(1+|\psi_j(\alpha, \beta)|)^{1+\rho}}\ d\beta\right)d\alpha+\right.\\
\non
&&\qquad\qquad\qquad\qquad+\left.\int^{\frac{t-r}2}_{0}\left(\int_{-\alpha}^{\alpha}\frac1{(1+\alpha)(t-r-\alpha)(1+|\psi_j(\alpha, \beta)|)^{1+\rho}}\ d\beta\right)d\alpha\right\}\\
\non
&\leq&\frac{C\{1+\log (1+t+r)\}M^{(i)}_{0, \rho}(H)}{(1+r)^{\frac12}(1+|t-r|)}.
\ea
This completes the proof of (\ref{l1}). \\
\section{$a$ $priori$ estimates}
In this section, we derive the $a$ $priori$ estimate (\ref{2.2}) assuming (\ref{2.1}).  For this purpose,  we introduce a notation. 
Let the assumptions of Lemma 2.1 be fulfilled and let $p(x, t)$ and  $q(x, t)$ be functions defined on a set $D\subset \mathbb{R}^2\times [0, T)$. Then we denote
$$
p(x, t)=O^*(q(x, t))\qquad\quad \mbox{in}\quad D,
$$
when there exist constants $K=K(B)>0$ and $\ve_0=\ve_0(J, B)>0$ such that, if (\ref{2.1}) holds for $0<\ve<\ve_0$, then 
$$
|p(x, t)|\leq K q(x, t) \qquad\mbox{for}\quad (x, t)\in D
$$
for the same $\ve$. We can easily show that if $p_1(x, t)=O^*(q(x, t))$ and $p_2(x, t)=O^*(q(x, t))$, then $p_1(x, t)+p_2(x, t)=O^*(q(x, t))$. Then our task to prove Lemma 2.1 is showing
\ba
 [\pa u(x, t)]_{k}=O^*(\ve)\quad\mbox{and}\quad \la u(x, t)\ra_{k+1}=O^*(\var)\qquad \quad \mbox{in}\quad \mathbb{R}^2\times [0, T_B).\label{4.1}
\ea
Also we will express constants determined independently of $J$  and $T$ by $K_n \ (n\in \mathbb{N})$ in the following argument. \\
\indent
Now we aim to show (\ref{4.1}). By (\ref{3.1}), we can write
\ba
u^i(x, t)=u_0^i(x, t)+L_{c_i}(F^i)(x, t),\label{4.2a}
\ea
where $u^i_0(x, t)$ is the solution to  (\ref{1.13}), (\ref{1.14}) and  satisfies for any nonnegative integer $p$, 
\ba
[[\pa u_0^i(t)]]_{p}+\la\la u_0^i(t)\ra\ra_{p+1}\leq C_0\ve\qquad \mbox{for}\quad  0\leq t<\infty,\label{4.2}
\ea
with some constant $C_0=C_0(f^i, g^i, p)>0$. 
Then, we have for a multi-index $a=(a_0, a_1, \cdots,a_4)$,
\ba
\Gamma^a L_{c_i}(F^i)=v_a^i+\sum_{|b|\leq |a|}C_{a, b}L_{c_i}(\Gamma^b F^i),\label{4.00}
\ea
with some constants $C_{a, b}$. Here, $v^i_a=v_a^i(x, t)$ is the solution to the Cauchy problem;
\ba
\non
\Box_i v_a^i=0\qquad \qquad\qquad\mbox{in}\quad \mathbb{R}^2\times [0, \infty),\\
\non
v_a^i(x, 0)=\ve^2\phi_a^i(x) ,\quad \pa_0v_a^i(x, 0)=\ve^2 \psi_a^i(x)\qquad \mbox{in}\quad \mathbb{R}^2,
\ea
with functions $\phi^i_a,\ \psi_a^i\in C_0^{\infty}(\mathbb{R}^2)$ determined by $(f^j,g^j)_{j=1, 2, \cdots,m}$ suitably.  
Indeed,  by the commutation relations of $\Gamma_{\alpha}$ and $\Box_i$ and by the definition of $L_{c_i}(F^i)$, we have
\ba
\non
\Box_i \Gamma_{\alpha} L_{c_i}(F^i)=\Gamma_{\alpha}\Box_iL_{c_i}(F^i)+2\delta_{\alpha 4}\Box_iL_{c_i}(F^i)=\Gamma_{\alpha}F^i+2\delta_{\alpha 4}F^i
\ea
and 
\ba
\Gamma_{\alpha} L_{c_i}(F^i)(x, 0)=0, \qquad\pa_0\Gamma_{\alpha} L_{c_i}(F^i)(x, 0)=\delta_{\alpha 0}F^i(x, 0).
\ea
Since $F^i$ is quadratic, we can denote $F^i(x, 0)=\ve^2 \psi^i(x)\in C^{\infty}_0(\mathbb{R}^2)$.  Hence we have
\ba
\Gamma_{\alpha}L_{c_i}(F^i)=v^i+L_{c_i}(\Gamma_{\alpha}F^i+2\delta_{\alpha 4}F^i)=v^i+L_{c_i}(\Gamma_{\alpha}F^i)+2\delta_{\alpha 4}L_{c_i}(F^i)
\ea
where $v^i=v^i(x, t)$ is the solution to the Cauchy problem;
\ba
\non
\Box_i v^i=0\qquad \qquad\qquad\mbox{in}\quad \mathbb{R}^2\times [0, \infty),\\
\non
v^i(x, 0)=0 ,\quad \pa_0v^i(x, 0)=\delta_{\alpha 0}\ve^2 \psi^i(x)\qquad \mbox{in}\quad \mathbb{R}^2.
\ea
This implies (\ref{4.00}) when $|a|=1$. Repeating the above argument, we can obtain (\ref{4.00}) for any $a$. \\
Note that, as with (\ref{4.2}), we have for any nonnegative integer $p$,
\ba
\sum_{|b|\leq p}[[\pa v_b^i(t)]]_{0}+\sum_{|c|\leq p+1}\la\la v_c^i(t)\ra\ra_{0}\leq C'_0\ve^2\qquad \mbox{for}\quad  0\leq t<\infty,\label{4.222}
\ea
with some constant $C'_0=C'_0(f^1,\cdots, f^m, g^1,\cdots, g^m,p)>0$. It follows from (\ref{4.2a}) and (\ref{4.00}) that 
\ba
\Gamma^a u^i(x, t)=\Gamma^a u^i_0(x, t)+v^i_a(x, t)+\sum_{|b|\leq |a|}C_{a, b}L_{c_i}(\Gamma^b F^i)(x, t)\label{4.01}
\ea
Therefore, our task for the proof of  (\ref{4.1}) is to show
\ba
\sum_{|b|\leq k} [\pa L_{c_i}(\Gamma^b F^i)(x, t)]_{0}+\sum_{|c|\leq k+1}\la L_{c_i}(\Gamma^c F^i)(x, t)\ra_{0}=O^*(\var)\qquad  \mbox{in}\quad \mathbb{R}^2\times [0, T_B).\label{4.3}
\ea 
We will show (\ref{4.3}) by dividing the area into 
some parts.\\
\indent
Firstly, we assume $\dis{0\leq t\leq 1/{\ve}}$. In this region, we can show the sharper estimates;
\ba
\begin{aligned}
&\sum_{|b|\leq k} [[\pa L_{c_i}(\Gamma^b F^i)(x, t)]]_{0}+\sum_{|c|\leq k+1}\la\la L_{c_i}(\Gamma^c F^i)(x, t)\ra\ra_{0}\\
=&\ O^*
(\var^{\frac54})\quad\qquad \mbox{in}\qquad\quad 
\mathbb{R}^2\times \big[0,  1/{\ve}\big].\label{4.3a}
\end{aligned}
\ea
For this purpose, we prepare two propositions with respect to the energy. \\
\begin{prop}
Let $u(x, t)\in C^2(\mathbb{R}^2\times [0, T))$ be a function satisfying $||u||_{2, T}<\infty$. Then, there exists a constant $C_3>0$ such that
\ba
|x|^{\frac12}|u(x, t)|\leq C_3||u(t)||_{2}
\label{4.4}
\ea
holds for $(x, t)\in \mathbb{R}^2\times [0, T)$.
\end{prop}
\begin{prop}
Let $\dis{u(x, t)=(u^1(x, t),u^2(x, t),\cdots,u^m(x, t))}$ $\dis{\in (C^{\infty}(\mathbb{R}^2\times [0, T)))^m}$ be the solution to (\ref{1.1}) and (\ref{1.2}) and also let $\ell$ be a positive integer.  Assume (\ref{1.7}). Then there exist constants $\delta>0$ and  $C_4=C_4(\ell)>0$ such that if $\dis{|\pa u|_{[\frac{\ell+1}2], T}<\delta}$ holds, 
then 
\ba
||\pa u(t)||_{\ell}\leq C_4||\pa u(0)||_{\ell}\exp\bigg(C_4\int_0^t|\pa u(s)|_{[\frac{\ell +1}2]}ds\bigg)
\label{4.5}
\ea
holds for $0\leq t<T$. 
\end{prop}
We omit the proof of the propositions. For the details of Proposition 4.1, see \cite{kla1}. On the other hand, we get Proposition 4.2 by the usual energy argument for the quasilinear wave equations. \\
By (\ref{2.1}) and $k\geq 21$, we have $|\pa u|_{[\frac{k+5}2], \frac1{\ve}}\leq |\pa u|_{k, \frac1{\ve}}\leq J\ve<\delta$ for $0<\ve <\ve_0$, if we take $\ve_0$ to be $J\ve_0<\delta$. Hence, by (\ref{2.1}) and (\ref{4.5}), we have
\ba
\non
||\pa u(t)||_{k+4}&=& O\bigg(C_4||\pa u(0)||_{k+4}\exp\bigg(C_4\int_0^{\frac1{\ve}}|\pa u(s)|_{[\frac{k+5}2]}ds\bigg)\bigg)\\
&=&O^*\bigg(\ve\exp\bigg(\int_0^{\frac1{\ve}}\frac{C_4J\ve}{\sqrt{1+s}}ds\bigg)\bigg)\label{4.6}\\
\non
&=& O^*\left(\ve\exp(4C_4J\ve^{\frac12})\right)\\
\non
&=& O^*(\ve) \qquad\qquad  \mbox{in}\quad \big[0,\ 1/{\ve}\big],
\ea
if we take $\ve_0$ to be $\ve_0\leq 1$ and $J^2\ve_0\leq 1$. Therefore, by (\ref{1.85}), (\ref{2.1}), (\ref{l0}), (\ref{4.4}) and (\ref{4.6}), we have for $|c|\leq k+1$
\ba
\non
&&\sum_{|c|\leq k+1}\la\la L_{c_i}(\Gamma^c F^i)(x, t)\ra\ra_{0}\\
\non
&=& O\left(\sum_{j=0}^m\sup_{(y, s)\in \atop \Lambda_j(\frac1{\ve})\cap D^i(x, t)}\{|y|^{\frac12}
(1+s+|y|)^{1+\mu}(1+||y|-c_js|)|F^i(y, s)|_{k+1}\}\right)\\
&=& O^*( (1+t)^{\frac9{16}+\mu}[\pa u(t)]_{[\frac{k+2}2]}||\pa u(t)||_{k+2})\label{4.7}\\
\non
&=&O^*(J\ve^{\frac{23}{16}-\mu})\\
\non
&=&O^*(\ve^{\frac54})\qquad\qquad \mbox{in}\qquad \mathbb{R}^2\times \big[0, \ 1/{\ve}\big],
\ea
if we take $\mu$ and $\ve_0$ to be $\mu<1/{16}$ and $J^8\ve_0\leq 1$. As for $\pa L_{c_i}(\Gamma^b F^i)$ with $|b|\leq k$, when $0\leq t\leq 1$, it follows from (\ref{1.85}) that $1+|x|+t$ and $1+||x|-c_it|$ are bounded in the support of the solution $u^i(x, t)$. Hence, by (\ref{4.00}), (\ref{4.222}) and (\ref{4.7}),  we find  
\ba
\non
\sum_{|b|\leq k}[[\pa L_{c_i}(\Gamma^b F^i)(x, t)]]_0&=&O\left(\sum_{|c|\leq k+1}([[v_c^i(x, t)]]_{0}+[[ L_{c_i}(\Gamma^{c}F^i(x, t)]]_{0})\right)\\
\non
&=&O\left(\sum_{|c|\leq k+1}(\la\la v_c^i(x, t)\ra\ra_{0}+\la\la L_{c_i}(\Gamma^{c}F^i(x, t)\ra\ra_{0})\right)\label{4.bb}\\
\non
&=&O^*(\ve^{\frac54})\qquad \qquad \mbox{in}\qquad \mathbb{R}^2\times [0, 1].
\ea
On the other hand, by (\ref{1.85}), (\ref{2.1}), (\ref{l1}) with $G=0$, (\ref{4.4}) and (\ref{4.6}), we have 
\ba
\non
&&\sum_{|b|\leq k}[[\nabla L_{c_i}(\Gamma^b F^i)(x, t)]]_{0}\\
\non
&=& O\left(\sum_{j=0}^m\sup_{(y, s)\in \atop \Lambda_j\cap D^i(x, t)}\{|y|^{\frac12}
(1+s+|y|)^{1+\nu}(1+||y|-c_js|)|F^i(y, s)|_{k+1}\}\right)\\
&=& O^*( (1+t)^{\frac9{16}+\nu}[\pa u(t)]_{[\frac{k+2}2]}||\pa u(t)||_{k+2})\label{4.10}\\
\non
&=&O^*(J\ve^{\frac{23}{16}-\nu})\\
\non
&=&O^*(\ve^{\frac54})\qquad\qquad \mbox{in}\qquad \mathbb{R}^2\times \big[0,\  1/{\ve}\big],
\ea
if we take $\nu$ and $\ve_0$ to be $\nu<1/{16}$ and $J^8\ve_0<1$. Therefore, when $1\leq t$, by (\ref{1.85}), (\ref{4.00}), (\ref{4.222}) , (\ref{4.7}), (\ref{4.10})  and the identity
\ba
\pa_0=-\frac{x_1}t\pa_1-\frac{x_2}t\pa_2+\frac1t S,\label{4.a2}
\ea
we have 
\ba
\begin{aligned}
&\ \sum_{|b|\leq k}[[\pa_0 L_{c_i}(\Gamma^b F^i)(x, t)]]_0\\
=&\ \ O\left(\sum_{|b|\leq k}\bigg([[\nabla  L_{c_i}(\Gamma^b F^i(x, t)]]_0
+\frac1{1+t}[[S L_{c_i}(\Gamma^b F^i)(x, t)]]_{0}\bigg)\right)\label{4.9}\\
=&\ \ O\left(\sum_{|b|\leq k}[[\nabla L_{c_i}(\Gamma^b F^i(x, t)]]_0
+\sum_{|c|\leq k+1}(\la\la v_c^i(x, t)\ra\ra_{0}+\la\la L_{c_i}(\Gamma^{c}F^i(x, t)\ra\ra_{0})
\right)\quad\\
=&\ \ O^*(\ve^{\frac54})\qquad \qquad \mbox{in}\qquad \mathbb{R}^2\times [1, 1/\ve].
\end{aligned}
\ea
Therefore we obtain (\ref{4.3a}).\\
\indent
Secondly, we assume $\dis{1/{\ve}\leq t\leq T_B}$. In this region, we need more precise energy estimate:
\begin{prop}
Let $u(x, t)=(u^1(x, t), u^2(x, t), \cdots, u^m(x, t))\in
(C^{\infty}( \mathbb{R}^2\times [0, T))^m$ be the solution to (\ref{1.1}) and (\ref{1.2}) under the same assumption in Theorem 1.1. Also let $\ell$ be a positive integer. Then there exist positive 
constants $C_5$ and $\delta$ such that if  
\ba
[[[\pa u]]]_{[\frac{\ell+1}2], T}+\la u\ra_{[\frac{\ell+1}2]+1,T}<\delta\label{delta}
\ea
holds, then 
\ba
||\pa u(t)||_{\ell}\leq C_5||\pa u(t_0)||_{\ell}\exp\left(C_5\int_{t_0}^t\frac{[[[\pa u(s)]]]_{[\frac{\ell+1}2]}+\la u(s)\ra_{[\frac{\ell+1}2]+1}}{1+s}\ ds\right)\label{4.11}
\ea
holds for $0\leq t_0 \leq t<T$. Here, we have set
\ba
\non
&&[[[v]]]_{p ,\tau}=\sup_{0\leq t< \tau}[[[v(t)]]]_{p},\quad [[[v(t)]]]_p=\sup_{x\in \mathbb{R}^2}[[[v(x, t)]]]_p,\\
\non
&&[[[v(x, t)]]]_p=\left\{\begin{array}{ll}[[v(x, t)]]_p& \mbox{when}\ \ |x|\leq t^{\frac{7}{8}},\\

[v(x, t)]_p& \mbox{when}\ \ |x|>t^{\frac7{8}}.\end{array}
\right.
\ea
\end{prop}
\vspace{3mm}
For the proof of (\ref{4.11}), see Proposition 4.1 in \cite{2006}, in there we used the $ghost$ $weight$ energy method.\\[2mm]
{\bf Remark:}\ The  difference of situations between Proposition 4.1 in \cite{2006} and Proposition 4.3 in this paper is the power of $1+||x|-c_it|$ in $[\pa u^i]_k$ which is the norm we  are going to estimate. To derive the $ghost$ $weight$ energy method, we need to suppose that $|\pa u^i(x, t)|_k$ is equivalent to $1+||x|-c_it|$ in the area where $|x|$ is small. Thus we need to define another norm $[[[\pa u]]]_k$ and assume  (\ref{delta}) in Proposition 4.3. \\[3mm]
\indent
In order to use (\ref{4.11}) 
with $\ell=k+9$ and $T=T_B$, we also show that
\ba
[[[\pa u(x, t)]]]_{[\frac{k+10}2]}+\la u(x, t)\ra_{[\frac{k+10}2]+1}=O^*
(\ve^{\frac12})\qquad \qquad \mbox{in}\qquad \mathbb{R}^2\times  [0, \ 
T_B)\label{4.13}
\ea
holds. 
 By (\ref{2.1}) and $k\geq 21$, we find that
\ba
[\pa u(x, t)]_{[\frac{k+10}2]}+\la u(x, t)\ra_{[\frac{k+10}2]+1}=O^*(J\ve) =O^*(\ve^{\frac12})\qquad \mbox{in}\qquad \mathbb{R}^2\times [0, T_B)
\label{4.13a}
\ea
if we take $\ve_0$ to be $J^2\ve_0\leq 1$. 
Furthermore, if we obtain 
\ba
\sum_{|b|\leq [\frac{k+10}2]}[[\nabla L_{c_i}(\Gamma^b F^i)(x, t)]]_{0}=O^*(\ve^{\frac12})\qquad\qquad \mbox{in}\qquad \mathbb{R}^2\times  \big[0, \ 
T_B\big),\label{4.12}
\ea  
then by (\ref{1.85}), (\ref{2.1}), (\ref{4.00}), (\ref{4.a2}) and (\ref{4.12}),  we find that
\ba
\begin{aligned}
&\ \ \sum_{|b|\leq [\frac{k+10}2]}[[\pa_0 L_{c_i}(\Gamma^b F^i)(x, t)]]_{0}\\
=&\ \ O\left(\sum_{|b|\leq [\frac{k+10}2]}\bigg(\frac{|x|}{t}[[\nabla L_{c_i}(\Gamma^b F^i)(x, t)]]_0+\frac1{t}[[SL_{c_i}(\Gamma^b F^i(x, t))]]_0\bigg)\right)\\
=&\ \ O\left(\sum_{|b|\leq [\frac{k+10}2]}[[\nabla L_{c_i}(\Gamma^b F^i)(x, t)]]_0+\right.\\
&\ \ \left.\qquad \qquad +\sum_{|c|\leq [\frac{k+10}2]+1}\bigg(\la\la v^i_c(x, t)\ra\ra_0+\frac{(1+|x|)^{\frac12}}{(1+|x|+t)^{\frac7{16}}}\la L_{c_i}(\Gamma^c F^i(x, t))\ra_0\bigg)\right)\\
=&\ \ O\left(\sum_{|b|\leq [\frac{k+10}2]}[[\nabla L_{c_i}(\Gamma^b F^i)(x, t)]]_0+\right.\\
&\ \ \left.\qquad\qquad\qquad+\sum_{|c|\leq [\frac{k+10}2]+1}(\la\la v^i_c(x, t)\ra\ra_0+\la L_{c_i}(\Gamma^c F^i(x, t))\ra_0)\right)\\
=&\ \ O^*(\ve^{\frac12}+J\ve)\\
=&\ \ O^*(\ve^{\frac12} )\qquad\qquad\mbox{in}\qquad \{(x, t)\ : \ |x|\leq t^{\frac78},\ 1/{\ve}\leq t<T_B\},\label{4.12a}
\end{aligned}
\ea
if we take $\ve_0$ to be $J^2\ve_0\leq 1$.  Hence, by (\ref{4.2}), (\ref{4.222}), (\ref{4.3a}),  (\ref{4.12}) and (\ref{4.12a}), we have (\ref{4.13}). \\
\indent
In order to prove (\ref{4.3}) and (\ref{4.12}), we show that for any positive integer $\ell\leq k+1$ and for any positive constant $\eta$
\ba
\non
&&\sum_{i=1}^m\sum_{|c|\leq \ell+1} \la\la L_{c_i}(\Gamma^c F^i)(x, t)\ra\ra_0\\
&=&O^*\big(\ve+J^2\ve^2(1+t)^{\eta}\sup_{0\leq s\leq t}||\pa u(s)||_{\ell+8}\big)\qquad\mbox{in}\qquad \mathbb{R}^2\times  \big[1/{\ve}, \ 
T_B\big)\label{4.14}
\ea
and
\ba
\non
&&\sum_{i=1}^m\sum_{|b|\leq \ell} [[\nabla L_{c_i}(\Gamma^b F^i)(x, t)]]_0\\
&=&O^*\big(\ve+J^2\sup_{(y,s)\in \mathbb{R}^2\times [0, t]}
|y|^{\frac12}|\pa u(y, s)|_{\ell+6}\big)\qquad\mbox{in}\qquad \mathbb{R}^2\times  \big[1/{\ve}, \ 
T_B\big).\label{4.15}
\ea
hold.  We will show (\ref{4.14}) and (\ref{4.15}) step by step. \\
\indent
At first, by (\ref{1.5}), (\ref{1.6}), (\ref{1.85}), (\ref{2.13}), (\ref{2.14}), 
(\ref{l0}), (\ref{l1}), (\ref{4.4}) and $\ve^2\log (1+t)\leq B$, we have for any $\mu_1>0$ and $\rho_1>0$, 
\ba
\begin{aligned}
&\ \ \sum_{i=1}^m\sum_{|c|\leq \ell+1} \la\la L_{c_i}(\Gamma^c F^i)(x, t)\ra\ra_0\\
=&\ \ O\bigg(\sum_{i=1}^{m}\sum_{|c|\leq \ell+1}
M_{1+\mu_1, 1}^{(i)}(\Gamma^c F^i)(x, t)
\bigg)\\
=&\ \ O\bigg(\sum_{i=1}^{m}\sum_{j=0}^m
\sup_{(y, s)\in D^i(x, t)\cap \Lambda_j}
|y|^{\frac12}z^{(j)}_{1+\mu_1, 1}(|y|, s)|F^i(y, s)|_{\ell+1}\bigg)\\
=&\ \ O\bigg(\sum_{j=0}^m \sup_{(y, s)\in \Lambda_j\atop 0\leq s\leq t}
\bigg\{|y|^{\frac12}z^{(j)}_{\mu_1,2}(|y|, s)\sum_{h=1}^m|\pa u^h(y, s)|_{[\frac{\ell+2}2]}|\pa u^h(y, s)|_{\ell+2}+\\
&\  \ 
+|y|^{\frac12}z^{(j)}_{\mu_1, 1}(|y|, s)\sum_{h=1}^m(|u^h(y, s)|_{[\frac{\ell+2}2]+1}|\pa u^h(y, s)|_{\ell+2}
+|\pa u^h(y, s)|_{[\frac{\ell+2}2]}|u^h(y, s)|_{\ell+3})+\\
&\ \ 
+|y|^{\frac12}z^{(j)}_{1+\mu_1, 1}(|y|, s)|\pa u(y, s)|_{[\frac{\ell+2}2]}^2|\pa u(y, s)|_{\ell+2}\bigg\}\bigg)\label{4.16}\\
=&\ \ O\bigg(([\pa u]_{[\frac{\ell+2}2],t}+\la u\ra_{[\frac{\ell+2}2]+1,t})\times\\
&\ \ \qquad\quad\times\sup_{(y, s)\in\mathbb{R}^2\times [0, t]}\{(1+s+|y|)^{-\frac7{16}+\mu_1}([[\pa u(y, s)]]_{\ell+2}+\la\la u(y,s)\ra\ra_{\ell+3})\}+\\
&\ \ \qquad +(1+t)^{\mu_1}[\pa u]_{[\frac{\ell+2}2],t}^2\sup_{(y,s)\in \mathbb{R}^2\times [0, t]}|y|^{\frac12}|\pa u(y, s)|_{\ell +2}
\bigg)\\
=&\ \ O^*\bigg(J\ve\sup_{(y, s)\in\mathbb{R}^2\times [0, t]}\{(1+s+|y|)^{-\frac7{16}+\mu_1}([[\pa u(y, s)]]_{\ell+2}+\la\la u(y,s)\ra\ra_{\ell+3})\}+\\
&\ \ \qquad +J^2\ve^2(1+t)^{\mu_1}\sup_{0\leq s\leq t}||\pa u(s)||_{\ell +4}
\bigg)\qquad\qquad\mbox{in}\qquad \mathbb{R}^2\times  \big[1/{\ve}, \ 
T_B\big)
\end{aligned}
\ea
and
\ba
\non
&&\sum_{i=1}^m\sum_{|b|\leq \ell} [[\nabla L_{c_i}(\Gamma^b F^i)(x, t)]]_0\\
\non
&=&O\bigg(\sum_{i=1}^{m}\sum_{|c|\leq \ell+1}
\bigg\{
M_{1+\mu_1, 1}^{(i)}(\Gamma^c N^i_2)
+(1+\log(1+t))M_{1, 1+\rho_1}^{(i)}(\Gamma^c (F^i-N_2^i) )
\bigg\}
\bigg)\\
\non
&=&O\bigg(\sum_{i=1}^{m}\sum_{j=0}^m
\sup_{(y, s)\in D^i(x, t)\cap \Lambda_j}
\{|y|^{\frac12}z^{(j)}_{1+\mu_1, 1}(|y|, s)|N_2^i(y, s)|_{\ell+1}+\\
\non
&&\qquad\qquad\qquad\qquad
+(1+\log (1+t))|y|^{\frac12}z^{(j)}_{1, 1+\rho_1}(|y|, s)|(F^i-N_2^i)(y, s)|_{\ell+1}\}\bigg)\\
\non
&=&O\bigg(\sum_{j=0}^m \sup_{(y, s)\in \Lambda_j\atop 0\leq s\leq t}
\bigg\{|y|^{\frac12}z^{(j)}_{\mu_1,2}(|y|, s)\sum_{h=1}^m|\pa u^h(y, s)|_{[\frac{\ell+2}2]}|\pa u^h(y, s)|_{\ell+2}+
\ea
\ba
\non
&&
+|y|^{\frac12}z^{(j)}_{\mu_1, 1}(|y|, s)\sum_{h=1}^m(|u^h(y, s)|_{[\frac{\ell+2}2]+1}|\pa u^h(y, s)|_{\ell+2}
+|\pa u^h(y, s)|_{[\frac{\ell+2}2]}|u^h(y, s)|_{\ell+3})\bigg\}+\\
\non
&&
+(1+\log (1+t))\sum_{j=0}^m \sup_{(y, s)\in \Lambda_j\atop 0\leq s\leq t}\{|y|^{\frac12}z^{(j)}_{1, 1+\rho_1}(|y|, s)|\pa u(y, s)|_{[\frac{\ell+2}2]}^2|\pa u(y, s)|_{\ell+2}\}\bigg)\\
&=&O\bigg(([\pa u]_{[\frac{\ell+2}2],t}+\la u\ra_{[\frac{\ell+2}2]+1,t})\times\label{4.17}\\
\non
&&\qquad\quad \times\sup_{(y, s)\in\mathbb{R}^2\times [0, t]}\{(1+s+|y|)^{-\frac7{16}+\mu_1}([[\pa u(y, s)]]_{\ell+2}+\la\la u(y,s)\ra\ra_{\ell+3})\}+\\
\non
&&\qquad +(1+\log(1+t))[\pa u]_{[\frac{\ell+2}2],t}^2\sup_{(y,s)\in \mathbb{R}^2\times [0, t]}|y|^{\frac12}|\pa u(y, s)|_{\ell +2}
\bigg)\\
\non
&=&O^*\bigg(J\ve\sup_{(y, s)\in\mathbb{R}^2\times [0, t]}\{(1+s+|y|)^{-\frac7{16}+\mu_1}([[\pa u(y, s)]]_{\ell+2}+\la\la u(y,s)\ra\ra_{\ell+3})\}+
\\
\non
&&\qquad +J^2\ve^2(1+\log(1+t))\sup_{(y,s)\in \mathbb{R}^2\times [0, t]}|y|^{\frac12}|\pa u(y, s)|_{\ell +2}
\bigg)\\
\non
&=&O^*\bigg(J\ve\sup_{(y, s)\in\mathbb{R}^2\times [0, t]}\{(1+s+|y|)^{-\frac7{16}+\mu_1}([[\pa u(y, s)]]_{\ell+2}+\la\la u(y,s)\ra\ra_{\ell+3})\}+\\
\non
&&\qquad +J^2\sup_{(y,s)\in \mathbb{R}^2\times [0, t]}|y|^{\frac12}|\pa u(y, s)|_{\ell +2}
\bigg)\qquad\mbox{in}\qquad \mathbb{R}^2\times  \big[1/{\ve},\ 
T_B\big),
\ea
where we have set 
\ba
\non
N_2^i=\sum_{j, \ell=1}^m\sum_{\alpha,\beta=0}^2 a_{\ell j}^{i, \alpha\beta\gamma}\pa_{\gamma}u^j\pa_{\alpha}\pa_{\beta}u^{\ell}
+\sum_{j,k=1}^m\sum_{\alpha,\beta=0}^2b_{jk}^{i, \alpha\beta}\pa_{\alpha}u^j\pa_{\beta}u^k.
\ea
Next, we estimate $(1+s+|y|)^{-\frac7{16}+\mu_1}([[\pa u(y, s)]]_{\ell+2}+\la\la u(y,s)\ra\ra_{\ell+3})$ for $(y, s)\in\mathbb{R}^2\times [0, t]$.  By the same 
manner as (\ref{4.16}), for any $\mu_2>0$, we obtain by  (\ref{1.5}), (\ref{1.6}), (\ref{1.85}), (\ref{2.13}), (\ref{2.14}), 
(\ref{l0}), (\ref{l1}), (\ref{4.2}), (\ref{4.222}), (\ref{4.01}) and (\ref{4.9})
\ba
\non
&& (1+s+|y|)^{-\frac7{16}+\mu_1}([[\pa u(y, s)]]_{\ell+2}+\la\la u(y,s)\ra\ra_{\ell+3})\\
\non
&=&O\bigg(\ve+\sum_{i=1}^m\sum_{|b|\leq \ell+2\atop |c|\leq \ell+3} (1+s+|y|)^{-\frac7{16}+\mu_1}\times\\
\non
&&\qquad\qquad\qquad\qquad\times([[\nabla L_{c_i}(\Gamma^b F^i)(y, s)]]_0+\la\la L_{c_i}(\Gamma^c F^i)(y, s)\ra\ra_0)\bigg)\\
\non
&=&O\bigg(\ve+\sum_{i=1}^{m}\sum_{|c|\leq \ell+3}\sup_{y\in\mathbb{R}^2}\{(1+s+|y|)^{-\frac7{16}+\mu_1}M_{1+\mu_2, 1}^{(i)}(\Gamma^c F^i)(y, s)\}
\bigg)\\
\non
&=&O\bigg(\ve+\sum_{i=1}^{m}\sum_{j=0}^m\sup_{y\in\mathbb{R}^2} \sup_{(\xi, \tau)\in D^i(y, s)\cap \Lambda_j}
\{|\xi|^{\frac12}z^{(j)}_{\frac9{16}+\mu_1+\mu_2, 1}(|\xi|, \tau)|F^i(\xi, \tau)|_{\ell+3}\bigg)
\ea
\ba
&=&O\bigg(\ve+([\pa u]_{[\frac{\ell+4}2],t}+\la u\ra_{[\frac{\ell+4}2]+1,t})\times\label{4.18}\\
\non
&&\qquad\qquad\quad\times \sup_{(\xi, \tau)\in\mathbb{R}^2\times [0, s]}\{(1+\tau+|\xi|)^{-\frac78
+\mu_1+\mu_2}([[\pa u(\xi, \tau)]]_{\ell+4}+\la\la u(\xi, \tau)\ra\ra_{\ell+5})\}+\\
\non
&&\qquad +[\pa u]_{[\frac{\ell+4}2],t}^2\sup_{(\xi, \tau)\in \mathbb{R}^2\times [0, s]}|\xi|^{\frac12}|\pa u(\xi, \tau)|_{\ell +4}
\bigg)\\
\non
&=&O^*\bigg(\ve+J\ve\sup_{(\xi, \tau)\in\mathbb{R}^2\times [0, s]}\{(1+\tau+|\xi|)^{-\frac78+\mu_1+\mu_2}([[\pa u(\xi, \tau)]]_{\ell+4}+\la\la u(\xi, \tau)\ra\ra_{\ell+5})\}+\\
\non
&&\qquad +J^2\ve^2\sup_{(\xi, \tau)\in \mathbb{R}^2\times [0, s]}|\xi|^{\frac12}|\pa u(\xi, \tau)|_{\ell +4}
\bigg)\qquad\mbox{in}\qquad \mathbb{R}^2\times  \big[1/{\ve}, \ 
T_B\big).
\ea
Moreover, by the same manner as (\ref{4.18}), for any $\mu_3>0$ we obtain  (\ref{1.5}), (\ref{1.6}), (\ref{1.85}), 
(\ref{l0}), (\ref{l1}), (\ref{4.2}), (\ref{4.222}), (\ref{4.01})  and (\ref{4.9})
\ba
\non
&& (1+\tau+|\xi|)^{-\frac7{8}+\mu_1+\mu_2}([[\pa u(\xi, \tau)]]_{\ell+4}+\la\la u(\xi, \tau)\ra\ra_{\ell+5})\\
\non
&=&O\bigg(\ve+\sum_{i=1}^m\sum_{|b|\leq \ell+4\atop |c|\leq \ell+5} (1+\tau+|\xi|)^{-\frac7{8}+\mu_1+\mu_2}\times\\
\non
&&\qquad\qquad\qquad\qquad\times([[\nabla L_{c_i}(\Gamma^b F^i)(\xi,  \tau)]]_0+\la\la L_{c_i}(\Gamma^c F^i)(\xi, \tau)\ra\ra_0)\bigg)\\
&=&O\bigg(\ve+\sum_{i=1}^{m}\sum_{|c|\leq \ell+5}\sup_{y\in\mathbb{R}^2}\{(1+\tau+|\xi|)^{-\frac7{8}+\mu_1+\mu_2}M_{1+\mu_3, 1}^{(i)}(\Gamma^c F^i)(\xi, \tau)\}
\bigg)\label{4.19}\\
\non
&=&O\bigg(\ve+\sum_{i=1}^{m}\sum_{j=0}^m\sup_{y\in\mathbb{R}^2} \sup_{(\zeta, \theta)\in D^i(\xi, \tau)\cap \Lambda_j}
\{|\zeta|^{\frac12}z^{(j)}_{\frac1{8}+\mu_1+\mu_2+\mu_3, 1}(|\zeta|, \theta)|F^i(\zeta, \theta)|_{\ell+5}\bigg)\\
\non
&=&O\bigg(\ve+[\pa u]_{[\frac{\ell+6}2],t} \sup_{(\zeta, \theta)\in\mathbb{R}^2\times [0, \tau]}\{|\zeta|^{\frac12}(1+\theta+|\zeta|)^{-\frac5{16}
+\mu_1+\mu_2+\mu_3}|\pa u(\zeta, \theta)|_{\ell+6}\}\bigg)\\
\non
&=&O^*\bigg(\ve+J\ve \sup_{(\zeta, \theta)\in\mathbb{R}^2\times [0, \tau]}\{|\zeta|^{\frac12}|\pa u(\zeta, \theta)|_{\ell+6}\}
\bigg)\qquad\mbox{in}\qquad \mathbb{R}^2\times  \big[1/{\ve}, \ 
T_B\big),
\ea
if we take $\mu_1,\mu_2,\mu_3$ to be $\mu_1+\mu_2+\mu_3<5/{16}$. Combining (\ref{4.16}), (\ref{4.17}), (\ref{4.18}) and (\ref{4.19}) and taking $\mu_1=\eta$, $\ve_0\leq 1$, $J\ve_0\leq 1$,
 we have (\ref{4.14}) and (\ref{4.15}).\\
\indent
Now we show (\ref{4.12}).  It follows from (\ref{2.1}), (\ref{4.15}) and $k\geq 21$  that
\ba
\non
&&\sum_{i=1}^m\sum_{|b|\leq [\frac{k+10}2]} [[\nabla L_{c_i}(\Gamma^b F^i)(x, t)]]_0\\
\non
&=&O^*\bigg(\ve+J^2B\sup_{(y,s)\in \mathbb{R}^2\times [0, t]}
|y|^{\frac12}|\pa u(y, s)|_{[\frac{k+10}2]+6}\bigg)\\
\non
&=&O^*(\ve+J^2[\pa u]_{k,t})\\
\non
&=&O^*((1+J^3)\ve)\\
\non
&=&O^*(\ve^{\frac12})\qquad\mbox{in}\qquad \mathbb{R}^2\times  \big[1/{\ve}, \ 
T_B\big),
\ea
if we take $\ve_0$ to be $J^6\ve_0\leq 1$.  This implies (\ref{4.12}) and therefore (\ref{4.13}). Furthermore, (\ref{4.13}) implies that there exists a positive constant $K_1$ such that 
\ba
[[[\pa u(s)]]]_{[\frac{k+10}2]}+\la u(s)\ra_{[\frac{k+10}2]+1}\leq K_1\ve^{\frac12}\label{4.13aa}
\ea
holds for $0<\ve <\ve_0$. Hence, by (\ref{4.11}) and (\ref{4.13aa}), we have
\ba
\non
||\pa u(t)||_{k+9}&=&O\bigg(C_5||\pa u(0)||_{k+9}\exp\bigg(C_5\int_0^t\frac{[[[\pa u(s)]]]_{[\frac{k+10}2]}+\la u(s)\ra_{[\frac{k+10}2]+1}}{1+s}\ ds
\bigg)\bigg)\\
\non
&=&O^*\bigg(\ve\exp\bigg(\ve^{\frac12}\int_0^t\frac{C_5K_1}{1+s}\ ds\bigg)\bigg)\\
&=&O^*\bigg(\ve\exp\bigg(C_5K_1\ve^{\frac12}\log (1+t)\bigg)\bigg)\label{4.20}\\
\non
&=&O^*\bigg(\ve(1+t)^{C_5K_1\ve^{\frac12}}\bigg)\qquad\qquad \mbox{in}\qquad \big[1/{\ve},\ T_B\big).
\ea
 Therefore , by 
(\ref{4.4}), (\ref{4.14}), (\ref{4.15}) and (\ref{4.20}), we obtain 
 \ba
 \non
 &&\sum_{i=1}^m\sum_{|c|\leq k+2}\la\la L_{c_i}(\Gamma^{c}F^i)(x, t)\ra\ra_{0}\\
\non
&=&O^*\bigg(\ve+J^2\ve^2(1+t)^{\eta}\sup_{0\leq s\leq t}||\pa u(s)||_{k+9}\bigg)\\
 &=&O^*\bigg(\ve+J^2\ve^3(1+t)^{\eta+C_5K_1\ve^{\frac12}}\bigg)\label{4.21}\\
 \non
 &=&O^*\big(\ve(1+t)^{\frac1{16}}\big)\qquad\qquad\mbox{in}\qquad \mathbb{R}^2\times  \big[1/{\ve}, \ 
T_B\big)
\ea
and
\ba
\non
&&\sum_{i=1}^m\sum_{|b|\leq k+1}[[\nabla L_{c_i}(\Gamma^b F^i)(x, t)]]_0\\
\non
&=&O^*\bigg(\ve+J^2\ve^2(1+\log(1+t))\sup_{0\leq s\leq t}||\pa u(s)||_{k+9}\bigg)\\
&=&O^*\bigg(\ve +J^2B\ve(1+t)^{C_5K_1\ve^{\frac12}}\bigg)\label{4.22}\\
\non
&=&O^*\big(\ve^{\frac{127}{128}}(1+t)^{\frac1{256}}\big)\qquad\qquad\mbox{in}\qquad \mathbb{R}^2\times  \big[1/{\ve}, \ 
T_B\big),
\ea
 if we choose $\eta$ and $\ve_0$ to be $0<\eta+C_5K_1\ve_0^{\frac12}<1/{16}$, $0<C_5K_1\ve^{\frac12}_0<1/256$ and $J^{256}\ve_0\leq 1$. Hence, 
 by (\ref{4.2}), (\ref{4.222}) and (\ref{4.21}), we have
 \ba
 \la\la u(x, t)\ra\ra_{k+2}=O^*(\ve(1+t)^{\frac1{16}})\qquad\qquad\mbox{in}\qquad \mathbb{R}^2\times  \big[1/{\ve}, \ 
T_B\big).\label{4.a10}
 \ea
Hence, by (\ref{4.a10}),we obtain 
\ba
\begin{aligned}
\la u(x, t)\ra_{k+2}&=O^*((1+t)^{-\frac1{16}}\la\la u(x, t)\ra\ra_{k+2})\\
&=O^*(\ve)\qquad\qquad\mbox{in}\qquad \mathbb{R}^2\times  \big[1/{\ve}, \ 
T_B\big)\label{4.a12}
\end{aligned}
\ea
and therefore by (\ref{4.2}), (\ref{4.222}), (\ref{4.01}), (\ref{4.9}), (\ref{4.22}) and (\ref{4.a12}), we have
\ba
 [\pa u(x, t)]_{k+1}=O^*\big(\ve^{\frac{127}{128}}(1+t)^{\frac1{256}}\big)
 \qquad\qquad\mbox{in}\qquad \mathbb{R}^2\times  \big[1/{\ve}, \ 
T_B\big).\label{4.a11}
\ea
Note that  (\ref{4.a12}) and (\ref{4.a11}) are stronger than we needed with respect to the order of derivatives. We will make use of the strength of  the estimates  below. \\
On the other hand, in order to estimate $\pa u$, we introduce a subset of $\mathbb{R}^2\times [0, T_B)$ by 
\ba
\non
\tilde{\Lambda}_i(T)=\big\{ (x, t)\ :\ ||x|-c_it|\leq t^{\frac14}, \ 1/{\ve}\leq t<T\big\}\qquad (i=1,2,\cdots, m),
\ea
and discuss by dividing the area $\mathbb{R}^2\times [0, T_B)$  into out-side and in-side of $\tilde{\Lambda}_i(T_B)$. We also introduce  notations;
\ba
\non
\tilde{\Lambda}^c_i(T)=\big\{ (x, t)\ :\ (x, t)\not\in \tilde{\Lambda}_i(T), \ \ 1/{\ve}\leq t<T\big\}
\ea
and 
\ba\
\non
\pa \tilde{\Lambda}_i(T)=\big\{ (x, t)\ :\ ||x|-c_it|=t^{\frac14}\ \ \mbox{when}\ \ 1/{\ve}<t<T\\
\non
\qquad\qquad\qquad\qquad\qquad \ \mbox{or}\ \ ||x|-c_it|<t^{\frac14}\ \ \mbox{when}\ \ t=1/{\ve}\ \big\}.
\ea
Then we find that
\ba
(1+t)^{\frac14}\leq C_6(1+||x|-c_i|)\qquad\quad \mbox{for}\qquad (x, t)\in \tilde{\Lambda}^c_i(T_B)\label{4.24}
\ea
holds for some constant $C_6>0$ and that
\ba
\tilde{\Lambda}_i(T_B)\subset \Lambda_i(T_B)\label{4.a24}
\ea
holds for sufficiently small $\ve>0$. 
Hence, it follows from (\ref{4.22}) and (\ref{4.24}) that 
\ba
\non
\sum_{|b|\leq k+1}[\nabla L_{c_i}(\Gamma^b F^i)(x, t)]_0
&=&O^*\big(\ve^{\frac{127}{128}}(1+t)^{\frac1{256}}(1+||x|-c_it|)^{-\frac1{16}}\big)\\
&=&O^*\big(\ve^{\frac{127}{128}}(1+t)^{-\frac3{256}}\big)\label{4.25}\\
\non
&=&O^*(\ve^{\frac{257}{256}})
\qquad\qquad\mbox{in}\qquad \tilde{\Lambda}_i^c(T_B),
\ea
for $i=1, 2, \cdots, m$. Hence, by (\ref{4.2}), (\ref{4.222}), (\ref{4.01}), (\ref{4.a2}), (\ref{4.a12}) and (\ref{4.25}), we obtain
\ba
[\pa u^i(x, t)]_{k+1}
&=&O^*(\ve)\label{4.b1}
\qquad\qquad\mbox{in}\qquad \tilde{\Lambda}_i^c(T_B),
\ea
for $i=1, 2, \cdots, m$. Especially, by (\ref{4.2a}), (\ref{4.2}), (\ref{4.00}), (\ref{4.222}), (\ref{4.a2}), (\ref{4.a12}) and (\ref{4.25}), we obtain
\ba
[\pa_0 u^i(x, t)-\ve\pa_0 u^i_0(x, t)]_0&=&O^*(\ve^{\frac{257}{256}})\qquad \qquad \mbox{on}\qquad \pa \tilde{\Lambda}_i(T_B)\ \label{4.a13}
\ea
and
\ba
[\pa^2_0 u^i(x, t)-\ve\pa^2_0 u^i_0(x, t)]_0&=&O^*(\ve^{\frac{257}{256}})\qquad \qquad 
\mbox{on}\qquad \pa \tilde{\Lambda}_i(T_B).\label{4.a14}
\ea
Now, the task left for us is to show 
\ba
[\pa u^i (x, t)]_k=O^*(\ve) \qquad\qquad\mbox{in}\qquad \tilde{\Lambda}_i(T_B),\label{4.26}
\ea
for $i=1, 2, \cdots, m$. We use the method of ordinary differential equation along the pseudo characteristic curves. Let $u(x, t)=(u^1(x, t), u^2(x, t),\cdots, u^m(x, t))$ be the solution to (\ref{1.1}) and (\ref{1.2}) and denote $x=r\omega, \ (r=|x|,\ \omega\in S^1)$. Then, for fixed $\lambda\in \mathbb{R}$ and $\omega\in S^1$, we define the $i$-th 
pseudo characteristic curve in $(r, t)$-plane by the solution $r=r^i(t; \lambda)$ of the Cauchy problem;
\ba
\frac{dr}{dt}&=&\kappa_i(r, t)\equiv c_i+\frac1{2c_i^3}\Theta_i(-c_i, \omega)(\pa_0 u^i(r\omega, t))^2
\qquad\quad t_0\leq t<T_B,\label{4.27}\\
\non
r(t_0)&=&c_it_0+\lambda,
\ea 
where $t_0=1/{\ve}$ when $|\lambda|<\ve^{-\frac14}$, and $t_0=\lambda^4$ when $|\lambda|\geq \ve^{-\frac14}$. Namely, the initial point $(r^i(t_0;\lambda)\omega, t_0)$ is on $\pa \tilde{\Lambda}_i(T_B)$ for each $\lambda\in\mathbb{R}$ and $\omega\in S^1$. Denote
\ba
\non
\mathcal{J}^i(\lambda;\omega)=\{(x, t)\ : \ x=r^i(t;\lambda)\omega,\ \ t_0\leq t<T_B\ \},
\ea
then we find that 
\ba
\non
\tilde{\Lambda}_i(T_B)=\dis{\bigcup_{\lambda\in\mathbb{R}, \ 
\omega\in S^1}}\mathcal{J}^i(\lambda; \omega)
\ea
holds for each $i=1, 2, \cdots, m$. For the details, see \cite{2005}. 
Now, we can transform the equation (\ref{1.1}) into an ordinary differential equation along the pseudo characteristic curve. For a vector valued function $v=(v^1, v^2, \cdots, v^m)$, set
\ba
\non
\mathcal{E}_iv=\Box_iv^i-\sum_{\ell=1}^{m}\sum_{\alpha,\beta=0}^2A_{\ell}^{i,\alpha\beta}(\pa u)\pa_{\alpha}\pa_{\beta}v^{\ell},
\ea
then we obtain an identity
\ba
\non
&&(\pa_0+\kappa_i\pa_r)\big(r^{\frac12}\pa_0v^i\big)\\
\non&=&\frac{r^{\frac12}}2\mathcal{E}_iv+\frac{r^{\frac12}}2(\pa_0+c_i\pa_r)^2v^i+
\frac{r^{\frac12}(\kappa_i-c_i)}{c_i}
(\pa_0+c_i\pa_r)\pa_0v^i+\\
&&+\frac{c_i^2}{2r^{\frac32}}\Omega^2v^i+\frac1{2r^{\frac12}}(\kappa_i-c_i)\pa_0v^i+\frac{c_i}{2r^{\frac12}}(\pa_0+c_i\pa_r)v^i\label{4.28}-\\
\non
&&-\frac{r^{\frac12}(\kappa_i-c_i)}{c_i}\pa^2_0 v^i+\frac{r^{\frac12}}2\sum_{\ell=1}^{m}\sum_{\alpha, \beta=0}^2A_{\ell}^{i, \alpha\beta}(\pa u)
\pa_{\alpha}\pa_{\beta}v^{\ell}.
\ea
Note that the  differential operator $\pa_0+\kappa_i\pa_r$ in the left hand side of (\ref{4.28}) means  the derivative  along $r=r^i(t;\lambda)$ in $(r, t)$-plane. Furthermore, by (\ref{2.6}), (\ref{2.7})  and the definition of $\tilde{\Lambda}_i(T)$, we have for any $\alpha,\beta=0,1,2$ and $i=1,2,\cdots, m$
\ba
\non
\frac1{r^{\frac12}}-\frac1{(c_it)^{\frac12}}&=&O\bigg(\frac{|r-c_it|}{(1+t)^{\frac32}}\bigg)=O\bigg(\frac1{(1+t)^{\frac54}}\bigg)\\
\non
(\pa_0+c_i\pa_r)v&=&O\bigg(\frac{1+|r-c_it|}{1+t}|\pa v|_0+\frac1{1+t}|v|_1\bigg)\\
\pa_{\alpha}v+\frac{\omega_{\alpha}}{c_i}\pa_0 v&=&O\bigg(\frac{1+|r-c_it|}{1+t}|\pa v|_0+\frac1{1+t}|v|_1\bigg)\label{4.29}\\
\non
(\pa_0+c_i\pa_r)^2v&=&O\bigg(\frac{(1+|r-c_it|)^2}{(1+t)^2}|\pa v|_1+\frac1{(1+t)^2}|v|_2\bigg)\\
\non
\pa_{\alpha}\pa_{\beta}v-\frac{\omega_{\alpha}\omega_{\beta}}{c_i^2}\pa_0^2v&=&
O\bigg(\frac{1+|r-c_it|}{1+t}|\pa v|_1+\frac1{(1+t)^2}|v|_2\bigg)
\ea
in $\tilde{\Lambda}_i(T_B)$. By (\ref{2.1}), (\ref{2.14}) and (\ref{4.29}),  we have
\ba
\non
&&
-\frac{\kappa_i-c_i}{c_i}\pa_0^2 v^i+\frac12\sum_{\ell=1}^{m}
\sum_{\alpha, \beta=0}^2A_{\ell}^{i, \alpha\beta}(\pa u)\pa_{\alpha}\pa_{\beta}v^{\ell}\\
\non
&=&-\frac1{2c_i^4}\Theta_i(-c_i, \omega)(\pa_0 u^i)^2\pa_0^2 v^i+\frac12\sum_{\alpha,\beta,\gamma,\delta=0}^2
c_{iii}^{i,\alpha\beta\gamma\delta}\pa_{\gamma}u^i\pa_{\delta}u^i\pa_{\alpha}\pa_{\beta}v^i+\\
&&+O^*\bigg(\frac{|r-c_it|}{1+t}|\pa u^i|_0|\pa v^i|_1+\frac1{1+t}|u^i|_1|\pa v^i|_1+\frac1{(1+t)^2}|u^i|_1|v^i|_2+\\
\non
&&+\sum_{j\ne i}(|\pa u^j|_0|\pa u|_0|\pa v|_1+|\pa u|_0^2|\pa v^j|_1)+|\pa u|^3_0|\pa v|_1\bigg)\label{4.30}\\
\non
&=&O^*\bigg(\frac{|r-c_it|}{1+t}|\pa u^i|_0|\pa v^i|_1+\frac1{1+t}|u^i|_1|\pa v^i|_1+\frac1{(1+t)^2}|u^i|_1|v^i|_2+\\
\non
&&\ \ \ +\sum_{j\ne i}(|\pa u^j|_0|\pa u|_0|\pa v|_1+|\pa u|_0^2|\pa v^j|_1)+|\pa u|_0^3|\pa v|_1
\bigg)\qquad \mbox{in}\qquad \tilde{\Lambda}_i(T_B),
\ea
if we take $\ve_0$ to be $J\ve_0\leq 1$. 
Therefore, it follows from (\ref{4.28}), (\ref{4.29})  and (\ref{4.30}) that
\ba
\begin{aligned}
&\ \ (\pa_0+\kappa_i\pa_r)\big(r^{\frac12}\pa_0 v^i\big)-\frac{r^{\frac12}}2\mathcal{E}_iv\\
=&\ \ O^*\bigg(\frac{1}{1+t}|\pa v^i|_1+\frac1{(1+t)^{\frac32}}|v^i|_2+\frac{1}{(1+t)^{\frac14}}|\pa  u^i|_0|\pa v^i|_1+\label{4.31}\\
&\ \ \ \ \ +\frac1{(1+t)^{\frac12}}|u^i|_1|\pa v^i|_1+r^{\frac12}\sum_{j\ne i}(|\pa u^j|_0|\pa u|_0|\pa v|_1+|\pa u|_0^2|\pa v^j|_1)+\\
&\ \ \ \ \ +r^{\frac12}|\pa u|^3_0|\pa v|_1\bigg)\qquad\mbox{in}\qquad \tilde{\Lambda}_i(T_B).
\end{aligned}
\ea
Now we show (\ref{4.26}) by induction with respect to $k$. Choose a point  $(x, t)\in \tilde{\Lambda}_i(T_B)$, then there exist $\lambda\in\mathbb{R}$ and  $\omega \in S^1$ such that  $x=r\omega$ and
$(r\omega, t)\in \mathcal{J}^i(\lambda;\omega)$. 
 At first, we show 
\ba
[\pa u^i(x, t)]_0=O^*(\ve)\qquad\qquad \mbox{in}\qquad \tilde{\Lambda}_i(T_B)\label{4.32}
\ea
for $i=1, 2, \cdots, m$. 
Setting $v=u$ in (\ref{4.31}), we have by (\ref{1.1}), (\ref{2.1}), (\ref{2.13}) and (\ref{2.15}),
\ba
\non
&&\frac{d}{ds}\big\{\big(r^i(s;\lambda)\big)^{\frac12}\pa_0u^i(r^i(s; \lambda)\omega, s)\big\}\\
\non
&=&O^*\bigg((r^i)^{\frac12}B^i(\pa u)+\frac{1}{1+s}|\pa u^i|_1+\frac1{(1+s)^{\frac32}}|u^i|_2+\frac{1}{(1+s)^{\frac14}}|\pa  u^i|_1^2+
\\
&&+\frac1{(1+s)^{\frac12}}|u^i|_1|\pa u^i|_1+(r^i)^{\frac12}\sum_{j\ne i}|\pa u^j|_1|\pa u|^2_1+(r^i)^{\frac12}|\pa u|^4_1\bigg)\label{4.33}\\
\non
&=&O^*\bigg(\frac{J\ve}{(1+s)^{\frac54}}\bigg)\qquad \mbox{in}\qquad [t_0, t],
\ea
if we take $\ve_0$ to be $J\ve_0\leq 1$. 
Integrating (\ref{4.32}) from $t_0$ to $t$, we have
\ba
\begin{aligned}
&r^{\frac12}\pa_0 u^i(r\omega, t)\\
=&\ \big(r^i(t_0;\lambda)\big)^{\frac12}\pa_0 u^i(r^i(t_0; \lambda)\omega, t_0)+O^*\bigg(\frac{J\ve}{(1+t_0)^{\frac14}}\bigg)\qquad\mbox{in}\qquad \tilde{\Lambda}_i(T_B),\label{4.34}\end{aligned}
\ea
which 
implies
\ba
r^{\frac12}\pa_0 u^i(r\omega, t)=O^*\big(\ve\big)\qquad\quad \mbox{in}\qquad \tilde{\Lambda}_i(T_B),\label{4.35}
\ea
if we take $\ve_0$ to be $J^4\ve_0\leq 1$. 
Moreover, integrating (\ref{4.27}) and using (\ref{4.35}), we have
\ba
\non
r-c_it&=&r^i(t_0; \lambda)-c_it_0+O^*\big(\ve^2\log (1+t)\big)\\
&=&r^i(t_0; \lambda)-c_it_0+O^*(B)\label{4.36}\\
\non
&=&O^*\big(1+|r^i(\lambda;t_0)-c_it_0|\big)\qquad\qquad\mbox{in}\qquad \tilde{\Lambda}_i(T_B).
\ea
Hence, by (\ref{4.25}), (\ref{4.34}) and (\ref{4.36}), we obtain
\ba
\non
&&[\pa_0 u^i(x, t)]_0\\
\non
&=&O\big((1+|r-c_it|)^{\frac{15}{16}}r^{\frac12}|\pa_0 u^i(r\omega, t)|\big)\\
\non
&=&O^*\bigg((1+|r^i(\lambda;t_0)-c_it_0|)^{\frac{15}{16}}\times\\
&&\qquad\qquad\times\bigg\{\big(r^i(t_0;\lambda)\big)^{\frac12}|\pa_0 u^i(r^i(\lambda;t_0)\omega, t_0)|+\frac{J\ve}{(1+t_0)^{\frac14}}\bigg\}\bigg)\label{4.37}\\
\non
&=&O^*\bigg(\ve+\frac{J\ve}{(1+t_0)^{\frac{1}{64}}}\bigg)\\
\non
&=&O^*(\ve)\qquad\qquad\mbox{in}\qquad\tilde{\Lambda}_i(T_B),
\ea
if we take $\ve_0$ to be $J^{64}\ve_0\leq 1$. It follows from  (\ref{2.1}), (\ref{4.29}) and (\ref{4.37})  that  (\ref{4.32}) holds. 
Note that (\ref{4.36}) implies that there exists a positive constant $C_7$ independent of $J$ such that
\ba
\frac1{C_7}(1+|r^i_0(t_0; \lambda)-c_it_0|)\leq 1+|r-c_it|\leq C_7(1+|r^i_0(t_0; \lambda)-c_it_0|)\label{4.a15}
\ea 
for $(r\omega, t)\in \mathcal{J}^i(\lambda;\omega)$. Therefore, by (\ref{1.19}), (\ref{4.a13}), (\ref{4.34}) and (\ref{4.a15}), we have
\ba
\non
&&(r^i(s;\lambda))^{\frac12}\pa_0 u^i(r^i(s; \lambda)\omega, s)\\
&=&-c_i\ve\pa_{\rho}\mathcal{F}(\lambda,\omega)+O^*\bigg(\frac{\ve^{\frac{257}{256}}}{(1+|r^i(t_0; \lambda)-c_i t_0|)^{\frac{15}{16}}}\bigg)\quad\qquad\mbox{in}\qquad [t_0, t].\label{4.a16}
\ea
Secondly, we will show
\ba
[\pa u^i(x, t)]_1=O^*(\ve)\qquad\qquad\mbox{in}\qquad \tilde{\Lambda}_i(T_B)\label{4.38}
\ea
for $i=1, 2, \cdots, m$. 
Set $v=\pa_0 u$ in (\ref{4.31}), then we have
\ba
(\pa_0+\kappa_i\pa_r)((r^i)^{\frac12}\pa^2_0 u^i)=\frac{(r^i)^{\frac12}}{2}\mathcal{E}_i\pa_0 u+O^*\bigg(\frac{J\ve}{(1+s)^{\frac54}}\bigg)
\qquad\qquad\mbox{in}\qquad [t_0, t].\label{4.39}
\ea
By (\ref{1.1}), (\ref{2.1}),  (\ref{2.13}), (\ref{2.14}), (\ref{2.15}), (\ref{4.29}) and (\ref{4.a16}), we have
\ba
\non
&&\mathcal{E}_i\pa_0 u\\
\non
&=&\sum_{\ell=1}^m\sum_{\alpha, \beta=0}^b\pa_0\big(A_{\ell}^{i,\alpha\beta}(\pa u)\big)\pa_{\alpha}\pa_{\beta}u^{\ell}+\pa_0\big(B^i(\pa u)\big)\\
\non
&=&\sum_{\alpha,\beta,\gamma,\delta=0}^2c_{iii}^{i, \alpha\beta\gamma\delta}(\pa_0\pa_{\gamma}u^i\pa_{\delta}u^i+
\pa_{\gamma}u^i\pa_0\pa_{\delta}u^i)\pa_{\alpha}\pa_{\beta}u^i+\\
\non
&&+O^*\bigg(\frac{1+|r^i-c_is|}{1+s}|\pa u^i|_1^2+\frac1{(1+s)^2}|\pa_0 u^i|_1|u^i|_2+\\
&&\qquad\qquad\qquad+\sum_{j\ne i}(|\pa u^j|_0|\pa u|_1^2+|\pa u|_0|\pa u^j|_1|\pa u|_1)+|\pa u|_0^2|\pa u|_1^2\bigg)\label{4.40}
\\
\non
&=&\frac{2\Theta_i(-c_i, \omega)}{c_i^4(r^i)^{\frac32}}((r^i)^{\frac12}\pa_0 u^i)((r^i)^{\frac12}\pa_0^2 u^i)^2+O^*\bigg(\frac{J^2\ve^2}{(1+s)^{\frac{39}{16}}}\bigg)\\
\non
&=&-\frac{2\ve \Theta_i(-c_i, \omega)\pa_{\rho}\mathcal{F}^i(\lambda, \omega)}{c_i^4(r^i)^{\frac12}(1+s)}((r^i)^{\frac12}\pa_0^2 u^i)^2+\\
\non
&&+O^*\bigg(\frac{J^2\ve^2}{(1+s)^{\frac{39}{16}}}+\frac{J^2\ve^{3+\frac1{256}}}{(1+s)^{\frac32}(1+|r^i(t_0; \lambda)-c_it_0|)^{\frac{15}{16}}}\bigg)\qquad\quad\mbox{in}\qquad [t_0, t].
\ea
Hence, we have
\ba
\non
&&(\pa_0+\kappa_i\pa_r)((r^i)^{\frac12}\pa^2_0 u^i)\\
&=&-\frac{\ve \Theta_i(-c_i, \omega)\pa_{\rho}\mathcal{F}^i(\lambda, \omega)}{c_i^4(1+s)}
((r^i)^{\frac12}\pa_0^2 u^i)^2+\\
\non
&&+O^*\bigg(\frac{J\ve}{(1+s)^{\frac{5}{4}}}+\frac{J^2\ve^2}{(1+s)^{\frac{31}{16}}}+\frac{J^2\ve^{3+\frac1{256}}}{(1+s)(1+|r^i(t_0; \lambda)-c_it_0|)^{\frac{15}{16}}}\bigg)
\qquad\quad\mbox{in}\qquad [t_0, t].\label{4.41}
\ea
Set 
\ba
W(s)=(r^i(s;\lambda))^{\frac12}\pa_0^2 u^i(r^i(s;\lambda)\omega, s),
\ea
then (\ref{1.19}), (\ref{4.a14}) and (\ref{4.41}) imply the Cauchy problem of the ordinary differential equation;
\ba
&&W'(s)=-\frac{\ve\Theta_i(-c_i, \omega)\pa_{\rho}\mathcal{F}^i(\lambda, \omega)}{c_i^4(1+s)}W(s)^2+Q(s),\qquad  t_0\leq s\leq t\ (<T_B),\ \ \ \label{4.42}\\
&&W(t_0)=(r^i(t_0, \lambda))^{\frac12}\pa_0^2 u^i(r^i(t_0, \lambda)\omega, t_0)=\ve c_i^2\pa_{\rho}^2\mathcal{F}^i(\lambda, \omega)+O^*\big(\ve^{\frac{257}{256}}\big),\ \ \ \ \label{4.43}
\ea
where
\ba
Q(s)=O^*\bigg(\frac{J\ve}{(1+s)^{\frac{5}{4}}}+\frac{J^2\ve^2}{(1+s)^{\frac{31}{16}}}+\frac{J^2\ve^{3+\frac1{256}}}{(1+s)(1+|r^i(t_0; \lambda)-c_it_0|)^{\frac{15}{16}}}\bigg).\label{4.44}
\ea
Note that 
\ba
\non
\int_{t_0}^t|Q(s)|\ ds&=&O^*\bigg( \frac{J\ve}{(1+t_0)^{\frac14}}+\frac{J^2\ve^2}{(1+t_0)^{\frac{15}{16}}}+\frac{J^3\ve^{3+\frac1{256}}\log (1+t)}{(1+|r^i(t_0; \lambda)-c_it_0|)^{\frac{15}{16}}}\bigg)\\
&=&O^*\bigg(\frac{J\ve^{\frac{17}{16}}}{(1+t_0)^{\frac{15}{64}}}+\frac{J^2\ve^{\frac{173}{64}}}{(1+t_0)^{\frac{15}{64}}}+\frac{J^3B\ve^{1+\frac1{256}}}{(1+|r^i(t_0; \lambda)-c_it_0|)^{\frac{15}{16}}}\bigg)\label{4.45}\\
\non
&=&O^*\bigg(\frac{\ve^{\frac{33}{32}}}{(1+t_0)^{\frac{15}{64}}}+\frac{\ve^{\frac{513}{512}}}{(1+|r^i(t_0; \lambda)-c_it_0|)^{\frac{15}{16}}}\bigg)
\qquad\qquad\mbox{in}\qquad [t_0,\ T_B),
\ea
if we choose $\ve_0$ to be $J^{3}\ve_0^{\frac1{512}}<1$.  Now we can show 
\ba
[\pa_0^2 u^i(x, t)]_0=O^*(\ve)\qquad\qquad\mbox{in}\qquad \tilde{\Lambda}_i(T_B),\label{4.46}
\ea
by using the following proposition.
\begin{prop}
Let $w(t)$ be the solution of the ordinary differential equation;
\ba
\non
w'(t)=\frac{\alpha}{1+t}w(t)^2+q(t)\qquad\qquad\mbox{for}\qquad T_0\leq t<T_1,
\ea
where $\alpha$ is a constant, $T_0$ and $T_1$  are positive constants and $q(t)$ is a continuous function in $[T_0, T_1)$.  Assume
\ba
\non
q_*=\int_{T_0}^{T_1}|q(t)|\ dt <\infty \qquad\mbox{and}\qquad 2\alpha q_*\{\log(1+T_1)-\log(1+T_0)\}<1.
\ea
Then, 
\ba
|w(t)|\leq \bigg(1+\frac1{1-\alpha (w(T_0)+q_*)\{\log(1+t)-\log (1+T_0)\} }\bigg)(|w(T_0)|+q_*)\label{4.ap}
\ea 
holds, as long as the right hand side of (\ref{4.ap}) is well-defined. 
\end{prop}
For the proof of Proposition 4.4, see the proof of Proposition 3.4 in \cite{2000}. \\
By (\ref{4.45}), we have $\dis{q_*=\int_{t_0}^{T_B}|Q(s)|\ ds=O^*(\ve^{\frac{513}{512}})<\infty}$ and 
\ba
\non
&&2\alpha q_*(\log (1+T_B)-\log (1+t_0))\\
\non
&=&-2\frac{\ve\Theta_i(-c_i, \omega)\pa_{\rho}\mathcal{F}^i(\lambda, \omega)}{c_i^4}q_*(\log(1+T_B)-\log (1+t_0))\\
\non
&\leq& K_{2}B\ve^{\frac1{512}}<1
\ea
for $0<\ve<\ve_0$, if we take $\ve_0$ to be $(K_{2}B)^{512}\ve_0<1$. Hence, it follows from (\ref{4.42}), (\ref{4.43}), (\ref{4.45}), (\ref{4.ap}) and $HB<1$ that
\ba
\non
&&|W(t)|\\
\non
&\leq& \left(1+\frac{1}{\dis{1-\alpha (W(t_0)+q_*)\{\log (1+t)-\log (1+t_0)\}}}\right)(|W(t_0)|+q_*)\\
&\leq& \left(1+\frac{1}{\dis{1-\bigg(\ve^2H 
+\alpha q_*\bigg)\{\log (1+t)-\log (1+t_0)\}}}\right)(|W(t_0)|+q_*)\label{4.47}\\
\non
&\leq& \left(1+\frac{1}{\dis{1-HB-K_{2}B\ve^{\frac1{512}}/2}}\right) (|W(t_0)|+q_*)\\
\non
&\leq& \left(1+\frac{2}{1-HB}\right) (|W(t_0)|+q_*)\qquad\qquad\quad t_0\leq t<T_B
\ea
holds for $0<\ve<\ve_0$, if we tale $\ve_0$ to be $\ve_0\leq \{(1-HB)/(K_{2}B)\}^{512}$.  Therefore, by (\ref{1.18}), (\ref{4.a15}), (\ref{4.43}) and (\ref{4.47}), we obtain 
\ba
\non
&&(1+|r-c_i t|)^{\frac{15}{16}}|W(t)|\\
\non
&=&
O^*\left(\bigg(1+\frac{2}{\dis{1-HB}}\bigg)C_7(1+|r^i(t_0;\lambda)-c_it_0|)^{\frac{15}{16}} (|W(t_0)|+q_*)\right)\\
\non
&=&O^*(\ve)\qquad\qquad\qquad \mbox{in}\qquad\tilde{\Lambda}_i(T_B).
\ea
This implies (\ref{4.46}). Moreover, by (\ref{2.1}), (\ref{4.29}) and (\ref{4.46}), we have
\ba
[\pa^2 u^i(x, t)]_0=O^*(\ve)\qquad\qquad\qquad\mbox{in}\qquad\tilde{\Lambda}_i(T_B).\label{4.kk}
\ea
\indent
Now we show (\ref{4.38}). Set $v=\Gamma_{p} u\ \ (p=3,4)$ in (\ref{4.31}), then we have
\ba
(\pa_0+\kappa_i\pa_r)((r^i)^{\frac12}\pa_0 \Gamma_{p} u^i)=\frac{(r^i)^{\frac12}}{2}\mathcal{E}_i\Gamma_{p} u+O^*\bigg(\frac{J\ve}{(1+s)^{\frac54}}\bigg)
\qquad\mbox{in}\qquad\tilde{\Lambda}_i(T_B).\label{4.48}
\ea
By (\ref{1.1}), (\ref{2.1}),  (\ref{2.13}), (\ref{2.14}), (\ref{2.15}), (\ref{4.29}), (\ref{4.32})  and (\ref{4.kk}), we have
\ba
\non
&&\mathcal{E}_i\Gamma_p u\\
\non
&=&\sum_{\ell=1}^m\sum_{\alpha, \beta=0}^2\big\{\Gamma_p \big(A_{\ell}^{i,\alpha\beta}(\pa u)\pa_{\alpha}\pa_{\beta}u^{\ell}\big)-A_{\ell}^{i,\alpha\beta}(\pa u)\pa_{\alpha}\pa_{\beta}\Gamma_p u^{\ell}\big\}-\\
\non
&&\ \ \  -[\Gamma_p, \Box_i]u^i+\Gamma_p \big(B^i(\pa u)\big)\\
&=&O^*\bigg(|\pa u|_0|\pa^2 u|_0|\pa_0 \Gamma_p u^i|_0+\frac{1+|r^i-c_is|}{1+s}|\pa u^i|_1^2+\frac1{(1+s)^2}|\pa_0 u^i|_1|u^i|_2+\label{4.49}\\
\non
&&\qquad +\sum_{j\ne i}(|\pa u^j|_0|\pa u|_1^2+|\pa u|_0|\pa u^j|_1|\pa u|_1)+|\pa u|_0^2|\pa u|_1^2\bigg)
\\
\non
&=&O^*\bigg(\frac{\ve^2}{(r^i)^{\frac12}(1+s)}|r^{\frac12}\pa_0\Gamma_p u^i|_0+\frac{J^2\ve^2}{(1+s)^{\frac{39}{16}}}\bigg)\qquad\quad\mbox{in}\qquad [t_0, t].
\ea
Hence, by setting
\ba
V_1(s)=(r^i(s;\lambda))^{\frac12}\pa_0\Gamma_p u^i(r^i(s;\lambda)\omega, s),
\ea
we have
\ba
V_1'(s)=O^*\bigg(\frac{\ve^2}{1+s}|V_1(s)|+\frac{J\ve}{(1+s)^{\frac54}}
\bigg)
\qquad \quad\mbox{in}\qquad [t_0, t],\label{4.50}
\ea
if we take  $\ve_0$\ to be  $J\ve_0\leq 1$. 
Thus, the Gronwall inequality implies
\ba
\non
&&|V_1(t)|\\
&=&O^*\bigg(\bigg\{ |V_1(t_0)|+\int_{t_0}^t  \bigg(\frac{J\ve}{(1+s)^{\frac54}}
\bigg) ds\bigg\}
\exp\bigg(\int_{t_0}^t\frac{K_{3}\ve^2}{1+s}ds\bigg)\bigg)\ \ \ \label{4.50a}\\
\non
&=&O^*\bigg(\bigg( |V_1(t_0)|+\frac{J\ve}{(1+t_0)^{\frac14}}
\bigg)
e^{K_{3}B}\bigg)\qquad\mbox{in}\qquad\tilde{\Lambda}_i(T_B).
\ea
Hence, by (\ref{4.b1}), (\ref{4.36}) and (\ref{4.50a}), we have
\ba
\non
&&(1+|r-c_it|)^{\frac{15}{16}}|V_1(t)|\\
\non
&=&O^*\bigg(C_7(1+|r^i(t_0, \lambda)-c_it_0|)^{\frac{15}{16}}\bigg( |V_1(t_0)|+\frac{J\ve}{(1+t_0)^{\frac14}}\bigg)
e^{K_{3}B}\bigg)\\
&=&O^*\bigg(\ve +\frac{J\ve}{(1+t_0)^{\frac1{64}}}\bigg)\label{4.51}\\
\non
&=&O^*(\ve)\qquad\qquad\mbox{in}\qquad\tilde{\Lambda}_i(T_B),
\ea
if we take $\ve_0$ to be $J^{64}\ve_0\leq 1$. Therefore, (\ref{4.a12}), (\ref{4.29}) and (\ref{4.51}) imply (\ref{4.38}). \\
\indent
Finally, for any integer $h$ so that $2\leq h\leq k$,  we show 
\ba
[\pa u^i (x, t)]_{h}=O^*(\ve)\qquad\qquad \qquad \mbox{in}\qquad \tilde{\Lambda}_i(T_B),\label{4.52}
\ea
for $i=1,\cdots,m$, under the assumption
\ba
[\pa u (x, t)]_{h-1}=O^*(\ve)\qquad\qquad \qquad \mbox{in}\qquad \mathbb{R}^2\times [1/\ve, T_B).
\label{4.53}
\ea
Set $v=\Gamma^a u$ with $|a|\leq h$ in (\ref{4.31}), then 
 (\ref{1.1}), (\ref{2.13}), (\ref{2.15}), (\ref{4.a12}), (\ref{4.a11}), (\ref{4.32}), (\ref{4.kk}) and (\ref{4.53})  
imply
\ba
\non
&&(\pa_0+\kappa_i\pa_r)((r^i)^{\frac12}\pa_0 \Gamma^a u^i)\\
\non
&=&\frac{(r^i)^{\frac12}}2\mathcal{E}_i \Gamma^a u+O^*\bigg(\frac{1}{1+s}|\pa u^i|_{h+1}+\frac1{(1+s)^{\frac32}}|u^i|_{h+2}+\\
&&\ \ +\frac{1}{(1+s)^{\frac14}}|\pa  u^i|_0|\pa u^i|_{h+1}+\frac1{(1+s)^{\frac12}}|u^i|_1|\pa u^i|_{h+1}+\label{4.54b}\\
\non
&&\ \ +(r^i)^{\frac12}\sum_{j\ne i}(|\pa u^j|_0|\pa u|_0|\pa u|_{h+1}+|\pa u|_0^2|\pa u^j|_{h+1})+(r^i)^{\frac12}|\pa u|^3_0|\pa u|_{h+1}\bigg)\\
\non
&=&\frac{(r^i)^{\frac12}}2\mathcal{E}_i \Gamma^a u+O^*\bigg(\frac{\ve}{(1+s)^{1+\frac{63}{256}}}\bigg)
\qquad\mbox{in}\qquad [t_0, t]
\ea
and
\ba
\non
&&\mathcal{E}_i \Gamma^a u\\
\non
&=&\sum_{\ell=1}^m\sum_{\alpha,\beta=0}^2\{\Gamma^a (A_{\ell}^{i, \alpha\beta}(\pa u)\pa_{\alpha}\pa_{\beta}u^{\ell})-
A_{\ell}^{i,\alpha\beta}(\pa u)\pa_{\alpha}\pa_{\beta}\Gamma^a u^{\ell}\}+\Gamma^a B^i(\pa u)+[\Box_i, \Gamma^a]u^i\\
\non
&=&O^*\bigg(|\pa u^i|_0|\pa^2 u^i|_0\sum_{|a|\leq h}|\pa_0 \Gamma^a u^i|+\frac{|r^i-c_is|}{1+s}(|\pa u^i|_h+|\pa u^i|_{h}^2)|\pa u^i|_{h+1}+\\
&&\ \ \ \ \ +\frac1{1+s}(|\pa u^i|_{h+1}+|\pa u^i|_{h+1}^2)|u^i|_{h+2}+|\pa u^i|_{h-1}^3+\label{4.54a}\\
\non
&&\ \ \ \ \ +\sum_{j\ne i}(|\pa u^j|_h|\pa u|_h|\pa u|_{h+1}+|\pa u|_h^2|\pa u^j|_{h+1})+|\pa u|^3_h|\pa u|_{h+1}\bigg)\\
\non
&=&O^*\bigg(\frac{\ve^2}{(r^i)^{\frac12}(1+s)}|(r^i)^{\frac12}\pa_0\Gamma^a u^i|+\frac{\ve^{1+\frac{127}{128}}}{(r^i)^{\frac12}(1+s)^{\frac{127}{256}}}+\\
\non
&&\ \ \ \ \ +\frac{\ve^3}{(r^i)^{\frac12}(1+s)(1+|r^i-c_is|)^{\frac{15}{16}}}\bigg)
\qquad\quad\mbox{in}\qquad [t_0, t].
\ea
Thus, by setting 
\ba
V_{h}(s)=\sum_{|a|\leq h}(r^i(s; \lambda))^{\frac12}\pa_0\Gamma^{a}u^i(r^i(s;\lambda)\omega, s)
\ea
and by (\ref{4.a15}), (\ref{4.54b}) and (\ref{4.54a}), we have
\ba
\non
&&V_h'(s)\\
\non
&=&O^*\bigg(\frac{\ve^2}{1+s}|V_h(s)|+\frac{\ve}{(1+s)^{1+\frac{63}{256}}}+\frac{\ve^3}{(1+s)(1+|r^i(t_0; \lambda)-c_it_0|)^{\frac{15}{16}}}\bigg)
\qquad\mbox{in}\quad [t_0, t].
\ea
The gronwall inequality implies
\ba
\non
&&|V_h(t)|\\
\non
&=&O^*\bigg(\bigg\{|V_h(t_0)|+\int_{t_0}^t\bigg(\frac{\ve}{(1+s)^{1+\frac{63}{256}}}+\frac{\ve^3}{(1+s)(1+|r^i(t_0;\lambda)-c_it_0|)^{\frac{15}{16}}}\bigg)ds\bigg\}\times\\
\non
&&\qquad\qquad\qquad\quad\times
\exp\bigg(\int_{t_0}^t \frac{K_{4}\ve^2}{1+s}ds\bigg)\bigg)\ \ \\
\non
&=&O^*\bigg(\bigg\{|V_h(t_0)|+\frac{\ve}{(1+t_0)^{\frac{63}{256}}}+\frac{B \ve}{(1+|r^i(t_0; \lambda)-c_it_0|)^{\frac{15}{16}}}\bigg\}
e^{ K_{4}B}\bigg)
\ea
Hence, by (\ref{4.a15}), we have
\ba
\non
&&(1+|r-c_i t|)^{\frac{15}{16}}|V_h(t)|\\
\non
&=&O^*\bigg(C_7(1+|r^i(t_0; \lambda)-c_it_0|)^{\frac{15}{16}}\times\\
&&
\qquad\times\bigg\{|V_h(t_0)|+\frac{\ve}{(1+t_0)^{\frac{63}{256}}}+\frac{B \ve}{(1+|r^i(t_0; \lambda)-c_it_0|)^{\frac{15}{16}}}\bigg\}
e^{ K_{4}B}\bigg)\\
\non
&=&O^*\bigg(\ve+\frac{\ve}{(1+t_0)^{\frac{3}{256}}}\bigg)\\
\non
&=&O^*(\ve)\qquad\qquad\qquad\mbox{in}\qquad \tilde{\Lambda}_i(T_B).\label{4.55}
\ea
It follows from (\ref{4.b1}) and (\ref{4.55}) that (\ref{4.52}) holds.

\end{document}